\documentclass[pdflatex,sn-basic,Numbered]{sn-jnl}

\usepackage{amsmath,amssymb,amsthm}
\usepackage{bm}
\usepackage{color}
\usepackage{enumitem}
\usepackage{subcaption}
\usepackage{graphicx}
\usepackage[section]{placeins}
\usepackage{float}
\usepackage{booktabs}
\usepackage{empheq}
\usepackage{tabularx}
\usepackage{geometry}
\theoremstyle{thmstyleone}
\newtheorem{theorem}{Theorem}[section]
\newtheorem{lemma}[theorem]{Lemma}

\newtheorem{corollary}[theorem]{Corollary}
\theoremstyle{thmstylethree}

\theoremstyle{thmstyletwo}
\newtheorem{remark}[theorem]{Remark}
\newtheorem{assumption}{Assumption}

\renewcommand{\vec}[1]{\bm{#1}}
\newcommand{\mat}[1]{\bm{#1}}

\begin{document}

\title[Optimal control near criticality]{Optimal control of symmetry-breaking dynamics near criticality}

\author*[1]{\fnm{Pearson W.} \sur{Miller}}\email{pwmiller@ucsd.edu}

\affil*[1]{\orgdiv{Department of Mathematics}, \orgname{University of California San Diego},
\orgaddress{\city{San Diego}, \state{CA}, \country{USA}}}

\abstract{Near criticality, dynamical systems reduce to normal forms, which unite nonlinear behavior across disciplines under well-known universality classes. We extend this analysis to optimal control, specifically to the optimal forcing of a system toward a target state near a cusp point, the generic codimension-two organizing center of bistability and of symmetry-breaking pitchfork bifurcations. Starting with optimality conditions obtained from the maximum principle, we conduct a formal asymptotic study of the resulting boundary-value problem and obtain the leading-order behavior of extremal solutions. Three asymptotic regimes are identified, corresponding to different scalings of the forcing strength. For intermediate and weak control, extremals lie close to a 2D extended center manifold outside thin boundary layers, and in the weak regime the unfolding of the cusp enters only through the reduced drift. Under this reduction, approximate closed-form solutions for the optimized trajectories are constructed using the method of matched asymptotics, with explicit front-interaction coefficients throughout the bistable part of the unfolding. We then derive the leading-order second variation in each regime. It is coercive for strong and intermediate control, while for weak control its sign is decided by a Jacobi equation in the reduced state, yielding explicit focal-point tests and showing that extremals with two or more turning points are never locally optimal. Our results are numerically validated against two classic systems: a bistable kinetic switch and a discrete reaction--diffusion system subject to a Turing instability.}



\maketitle

\section{Introduction}
\label{sec:introduction}

A bifurcation strips a system down to a normal form. Near criticality, after transients, the dynamics of a high-dimensional model are captured, up to small corrections, by a scalar amplitude equation, and system-specific detail survives only through a few coefficients. This paper extends this well-known reduction to an optimally forced trajectory near bifurcation. We focus on the control of a system near a cusp point, the generic codimension-two organizing center of bistability. Its unfolding contains both the pitchfork, the most basic of symmetry-breaking phenomena, and its imperfect, asymmetric perturbations. We investigate the optimal external forcings required to steer the system toward a target state over a finite horizon, under quadratic tracking penalties. Our goal spans two questions. First, under what conditions does the optimal control problem inherit the center-manifold behavior of the unforced system? Second, what is the form of optimal trajectories on this manifold?

\begin{figure}[!t]
    \centering
    \includegraphics[width=\textwidth]{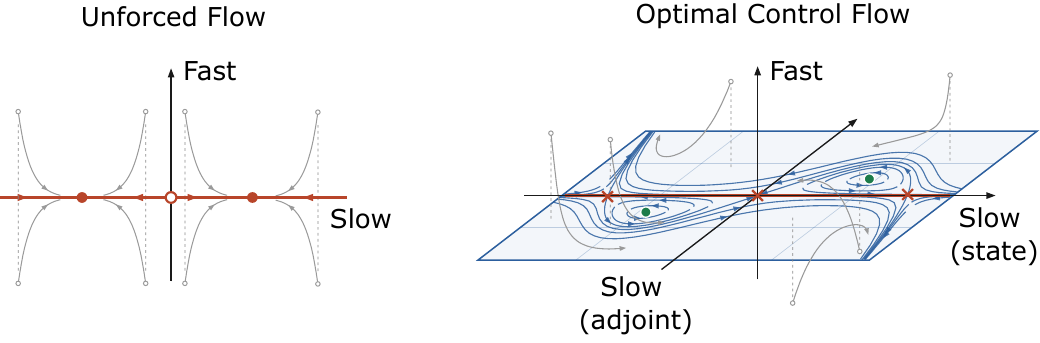}
    \caption{Schematic of the extended center-manifold problem for our optimal control problem, shown for the symmetric (pitchfork) case. Near the organizing center, a system possesses a 1D slow manifold (left) on which the dynamics are well approximated by a scalar amplitude equation. The optimal control problem is posed on the full system, but outside thin initial and terminal layers the optimal trajectory lies close to a 2D extended center manifold (right) that includes the original slow mode and its associated costate.}
    \label{fig:optimal_control_pitchfork}
\end{figure}

Bifurcation control is a mature subject, but its aims have historically been the reverse of ours. In most engineering scenarios, bifurcations relate to undesirable instabilities or performance limits, for instance unwanted vibrational instabilities in airfoils \cite{lee1999nonlinear}, and design or control is used to avoid or remove them \cite{monnigmann2002normal}. By contrast, in natural systems, especially in biology \cite{mora2011biological}, it is not unusual for tightly regulated dynamics to occur in a near-bifurcation regime. The examples are far-flung, though the most famous instance emerges from epigenetics, where cell phenotype transitions have long been associated with local bifurcations and the associated potential flows \cite{waddington1957strategy, ferrell2012bistability,wang2011quantifying}. Optimal control is increasingly being employed not just as a tool for design, but as a principled way to examine performance bounds and cost tradeoffs in natural systems \cite{pezzotta2023optimal, sinigaglia2024optimal, krishnan2026hamiltonian}. This growing interest establishes a need for better analytical characterization of the solutions of nonlinear optimal control problems.

We specify our problem as follows. Consider the controlled dynamical system
\begin{equation}
\frac{d c^i}{d t} = r^i(c^j, u^\alpha, h),  \quad t \in [0, T],
\label{eq:state_system}
\end{equation}
where $c^i(t) \in \mathbb{R}^N$ is the state, $u^\alpha(t) \in \mathbb{R}^M$ is the control input, and $h = (h_1, h_2) \in \mathbb{R}^2$ are parameters. The right-hand side $r^i: \mathbb{R}^N \times \mathbb{R}^M \times \mathbb{R}^2 \to \mathbb{R}^N$ is a sufficiently smooth nonlinear vector field, chosen so that the unforced system has a cusp point at $(c, h) = (0, 0)$ (made precise in Sec.~\ref{sec:preliminaries}).
We define the quadratic cost functional
\begin{equation}
J[u^\alpha] = \frac{\chi}{2} \|c(T) - c^{*}\|^2 + \frac{1}{2}\int_0^T \left[ \zeta \|c(t) - c^{*}\|^2 + \|u(t)\|^2 \right] dt,
\label{eq:cost_functional}
\end{equation}
where $c^{*i} \in \mathbb{R}^N$ is a prescribed target state, $\chi \ge 0$ and $\zeta > 0$ are the cost penalties, and $T > 0$ is the time horizon. Our goal is to characterize the asymptotic behavior of optimal controls
\begin{equation}
u^\alpha_{\text{opt}} \in \arg\min_{u^\alpha \in \mathcal{U}_{\mathrm{ad}}} J[u^\alpha],
\qquad
\mathcal{U}_{\mathrm{ad}} = \big\{ u^\alpha \in L^2([0,T],\mathbb{R}^M) : \text{the solution of \eqref{eq:state_system} exists on } [0,T] \big\},
\label{eq:optimization_problem}
\end{equation}
subject to the dynamics in Eq.~\eqref{eq:state_system} with initial condition $c^i(0) = c^i_0$. Minimizers need not be unique, since the problem is nonconvex. We analyze extremals, i.e.\ solutions of the first-order necessary conditions of Pontryagin's maximum principle.

This paper is concerned with what can be obtained about this problem via perturbative techniques. An important antecedent to the current work is that of Abed and Fu, who used bifurcation formulae to construct local feedback laws that stabilize an equilibrium at criticality and the bifurcated equilibria emanating from it \cite{abed1986local,abed1987local}. Kang and Krener subsequently placed this line of work in a systematic normal-form framework. They classified nonlinear control systems with uncontrollable linearization and characterized the resulting control bifurcations, degeneracies occurring where linear stabilizability is lost \cite{kang1998bifurcation1,krener2004control}. Hamzi, Kang, and Krener extended this to a controlled center-manifold reduction \cite{hamzi2005controlled}.

Perturbative constructions of nonlinear optimal control likewise have a long history. Al'brekht \cite{al1961optimal} and Lukes \cite{lukes1969optimal} solved the Hamilton--Jacobi--Bellman equation by expanding the value function in a power series about a stabilizable equilibrium. This construction has been developed extensively in a variety of ways, both numerically and analytically \cite{navasca2007patchy, sakamoto2008analytical, krener2025brekht}. These expansions are not a model reduction: they retain the full state space and do not isolate the intrinsic slow dynamics of the uncontrolled system. Here we seek instead a reduction that preserves those intrinsic dynamics, namely the reduced drift $\sigma h_1 + \mu h_2 A + \nu A^3$ on the center manifold (coefficients defined in Sec.~\ref{sec:preliminaries}).

The approach employed below is a simultaneous scaling of the state, the parameters and the control strength. The quasi-homogeneous cusp scaling $c \sim \varepsilon^{1/2}$, $h_2 \sim \varepsilon$, $h_1 \sim \varepsilon^{3/2}$ is imposed on the state and parameters. A matching scaling is imposed on the tracking penalties, weakening the overall control strength. This rescaling is opposite to the one typically studied, namely the cheap-control problem, for which well-developed theory exists \cite{kokotovic1999singular}.

The approach yields a three-regime structure determined by the scaling of $(\chi,\zeta)$ with $\varepsilon$, which we call the strong, intermediate, and weak control cases (Table~\ref{scal_tab}). In each regime the terminal and running penalties are scaled together so that $\chi^2/\zeta = O(1)$, and the horizon $T$ is taken commensurate with the timescale of that regime listed in Table~\ref{scal_tab}. We state the three regimes as a single reduction theorem for the Pontryagin system.
Under \emph{strong control}, $u = O(\varepsilon^{1/2})$, the optimal forcing generically excites the noncritical modes at leading order, and no general center-manifold reduction applies. Under \emph{intermediate control}, $u = O(\varepsilon)$, the center manifold exists, but the forcing remains strong enough to dominate the intrinsic drift. Under \emph{weak control}, $u = O(\varepsilon^{3/2})$, the control and the nonlinear drift are balanced. The result is a reduced nonlinear boundary-value problem on the center manifold (Fig.~\ref{fig:optimal_control_pitchfork}), which the optimal trajectory follows outside thin initial and terminal layers, and in which the unfolding parameters enter only through the reduced drift.

\begin{table}[htbp]
    \centering
    \begin{tabularx}{\textwidth}{llllX}
        \toprule
        Regime & $(\chi,\zeta)$ scaling & Control $u^\alpha$ & Timescale & Reduced system \\
        \midrule
        Strong control   & $(O(1), O(1))$       & $O(\varepsilon^{1/2})$  & $T \sim 1$              & LQR tracking for the linearization on $\mathbb{R}^N$ \\
        Intermediate control & $(O(\varepsilon^{1/2}), O(\varepsilon))$ & $O(\varepsilon)$    & $T \sim \varepsilon^{-1/2}$ & Driftless 1D LQR\\
        Weak control     & $(O(\varepsilon), O(\varepsilon^{2}))$ & $O(\varepsilon^{3/2})$ & $T \sim \varepsilon^{-1}$ & 1D nonlinear optimal control problem (2D Hamiltonian boundary-value problem)\\
        \bottomrule
    \end{tabularx}
     \caption{Summary of the three asymptotic scaling regimes. Each regime is defined by the joint scaling of the state tracking penalties $\chi$ and $\zeta$ with $\varepsilon$, and is characterized by a distinct rescaling of time. The corresponding leading-order optimal feedback law is qualitatively different in each case.}
     \label{scal_tab}
\end{table}

The remainder of this paper is organized as follows. Sec.~\ref{sec:preliminaries} states the cusp genericity conditions, fixes the parameter coordinates and the distinguished scaling, and derives the first-order optimality conditions. Sec.~\ref{sec:asymptotics} states the three regimes as cases of a single reduction theorem for the optimality system and derives each case explicitly, and Sec.~\ref{sec:second_order} gives the leading-order second variation. Sec.~\ref{sec:nonlinear_analysis} analyzes the nonlinear weak regime: a change of variables brings the reduced Hamiltonian to natural form, we carry out a bifurcation analysis on the pitchfork ray, and a matched asymptotic analysis yields closed-form approximations for the long-timescale behavior of extremals throughout the bistable part of the unfolding. Sec.~\ref{sec:conclusion} concludes. Appendix~\ref{app:expansions} collects the explicit coefficient formulas and expansions, and Appendix~\ref{app:numerics} contains details of the numerical validation.

\section{Preliminaries}
\label{sec:preliminaries}

In this section, we collect the elementary properties of our problem, list assumptions, and establish basic notation for use below. The unforced dynamics are organized by a codimension-two cusp point, and we follow the standard treatment of the cusp bifurcation and center-manifold reduction \cite{guckenheimer1983nonlinear, kuznetsov2023elements, kuznetsov1999numerical, mielke1991hamiltonian}.

\subsection{The organizing center}
\label{sec:organizing_center}

No particular structure is imposed on how $h$ enters $r^i$ in Eq.~\eqref{eq:state_system}. We assume that the unforced system has an equilibrium at the parameter value where the bifurcation takes place, and translate coordinates so that both sit at the origin:
\begin{equation}
r^i(0, 0, 0) = 0.
\label{eq:organizer}
\end{equation}
We refer to $(c^i, u^\alpha, h) = (0,0,0)$ as the organizing center.

\begin{remark}
Throughout, superscript Latin indices $i,j,k,l,m,n \in \{1,\ldots,N\}$ label components of $N$-vectors (state, costate, eigenvectors, etc.), and superscript Greek indices $\alpha,\beta,\gamma \in \{1,\ldots,M\}$ label components of $M$-vectors (control). Subscripts denote partial derivatives: $r^i{}_j \equiv \partial r^i/\partial c^j$, $r^i{}_{jk} \equiv \partial^2 r^i/\partial c^j\partial c^k$, and so on, with Greek subscripts for control derivatives ($r^i{}_\alpha \equiv \partial r^i/\partial u^\alpha$) and the subscripts $h_1, h_2$ for parameter derivatives ($r^i{}_{h_a} \equiv \partial r^i/\partial h_a$, $a \in \{1,2\}$). Unless explicitly noted, all such terms are evaluated at the organizing center. We follow the Einstein summation convention for Latin and Greek indices. The parameter label $a$ and the mode labels introduced below are never summed implicitly.
\end{remark}

\begin{assumption}[Simple critical eigenvalue]
\label{assump:bifurcation}
At the organizing center, the Jacobian $r^i{}_j$ has an algebraically simple zero eigenvalue, denoted $\lambda_c = 0$, and all its other eigenvalues have strictly negative real parts.
\end{assumption}

\begin{assumption}[Diagonalizability]
\label{assump:diag}
$r^i{}_j$ can be diagonalized over $\mathbb{C}$.
\end{assumption}
Let $v_c^i$ and $w_c^i$ denote real right and left eigenvectors for $\lambda_c = 0$, so that $r^i{}_j v_c^j = 0$ and $w_c^i r^i{}_j = 0$. For the stable eigenvalues $\lambda_s$ ($s = 2, \ldots, N$), let $v_s^i$ be right eigenvectors, and define $w_s^i$ by $\overline{w_s}^{\,i} r^i{}_j = \lambda_s \overline{w_s}^{\,j}$, where the overbar denotes complex conjugation; $\overline{w_s}^{\,i}$ then extracts the coefficient of mode $s$. With mode labels $\ell, \ell' \in \{c, 2, \ldots, N\}$, the eigenvectors are normalized by $\|v_\ell\| = 1$ and $\overline{w_\ell}^{\,i} v_{\ell'}^i = \delta_{\ell\ell'}$.

Three projected coefficients determine the leading-order dynamics on the center manifold. For $a \in \{1,2\}$, the coefficient $\sigma_a \equiv w_c^m r^m{}_{h_a}$ is the direct forcing of the critical mode by the parameter $h_a$, and $\mu_a$ is the rate at which $h_a$ changes the linear coefficient of the reduced vector field. The coefficient $\nu$ is the cubic coefficient of the reduced vector field. Both $\mu_a$ and $\nu$ include contributions from stable modes slaved to the critical mode; they are the standard center-manifold expressions \cite{kuznetsov1999numerical, kuznetsov2023elements}, given explicitly in Appendix~\ref{app:expansions}.

\begin{assumption}[Cusp point]
\label{assump:cusp}
The quadratic coefficient of the critical mode vanishes and the cubic coefficient does not:
\begin{equation}
  w^i_c\, r^i{}_{jk}\, v^j_c v^k_c = 0,
  \qquad
  \nu \neq 0.
  \label{eq:cusp_conditions}
\end{equation}
\end{assumption}

\begin{assumption}[Transversal unfolding]
\label{assump:unfolding}
The parameters unfold the cusp transversally:
\begin{equation}
  \det \mathsf{M} \neq 0,
  \qquad
  \mathsf{M} \equiv
  \begin{pmatrix}
    \sigma_1 & \sigma_2 \\
    \mu_1 & \mu_2
  \end{pmatrix}.
  \label{eq:transversality}
\end{equation}
\end{assumption}

Assumptions~\ref{assump:bifurcation}, \ref{assump:cusp} and \ref{assump:unfolding} are the standard genericity conditions for a cusp bifurcation of equilibria \cite[Secs.~8.2 and 8.7]{kuznetsov2023elements}. The first identity in \eqref{eq:cusp_conditions} is the codimension-two degeneracy: it makes the organizing center the point at which two fold (saddle-node) curves meet. The condition $\nu \neq 0$ is the nondegeneracy condition, which makes the cubic term the leading nonlinearity on the center manifold. Condition \eqref{eq:transversality} requires the map from $(h_1, h_2)$ to the constant and linear coefficients of the reduced vector field to be a local diffeomorphism at $h = 0$, so that the parameters generate both coefficients independently.
Under these conditions, the restriction of the unforced system to its center manifold is locally topologically equivalent to
\begin{equation*}
\dot{x} = \beta_1 + \beta_2\, x + \operatorname{sgn}(\nu)\, x^3,
\end{equation*}
with $(\beta_1,\beta_2)$ a smooth, invertible function of $(h_1,h_2)$ near the origin \cite{kuznetsov2023elements}.

The parameter-dependent quadratic coefficient of the critical mode, generically of size $O(|h|)$, does not appear in Assumption~\ref{assump:unfolding}. It can be removed by an $h$-dependent shift of the critical coordinate, and that shift alters the constant and linear coefficients only at second order in $h$. Consequently, the Jacobian of the parameter map at $h=0$, and with it the transversality condition, is unaffected.

It is convenient to normalize the parameter coordinates. Assumption~\ref{assump:unfolding} implies $(\sigma_1, \sigma_2) \neq (0,0)$. Labeling the parameters so that $\sigma_1 \neq 0$, consider the linear change of parameter coordinates
\begin{equation}
    h_1^{\prime} = h_1 + \frac{\sigma_2}{\sigma_1} h_2, \qquad h_2^{\prime} = h_2 .
    \label{eq:param_change}
\end{equation}
Inverting, $h_1 = h_1^{\prime} - (\sigma_2/\sigma_1)h_2^{\prime}$ and $h_2 = h_2^{\prime}$, so by the chain rule the parameter derivatives transform as
\begin{equation*}
    \frac{\partial}{\partial h_1^{\prime}} = \frac{\partial}{\partial h_1},
    \qquad
    \frac{\partial}{\partial h_2^{\prime}} = \frac{\partial}{\partial h_2} - \frac{\sigma_2}{\sigma_1}\frac{\partial}{\partial h_1}.
\end{equation*}
Since $\sigma_a$ and $\mu_a$ are linear in $\partial/\partial h_a$, the projected coefficients in the new coordinates are
\begin{equation*}
    \sigma^{\prime}_1 = \sigma_1, \qquad
    \sigma^{\prime}_2 = \sigma_2 - \frac{\sigma_2}{\sigma_1}\sigma_1 = 0, \qquad
    \mu^{\prime}_1 = \mu_1, \qquad
    \mu^{\prime}_2 = \mu_2 - \frac{\sigma_2}{\sigma_1}\mu_1 = \frac{\det \mathsf{M}}{\sigma_1}.
\end{equation*}
The transformation \eqref{eq:param_change} has unit determinant, so $\det \mathsf{M}^{\prime} = \sigma^{\prime}_1 \mu^{\prime}_2 = \det \mathsf{M}$. In the new coordinates, $h_2^{\prime}$ is the parameter direction that leaves the constant term unchanged at linear order, and $h_1^{\prime}$ carries all of it.

We adopt these coordinates from here on and drop the primes. This involves no loss of generality: it is a relabeling of parameter space, and the dependence of $r^i$ on $h$ remains otherwise arbitrary. Writing
\begin{equation}
    \sigma \equiv \sigma_1 = w_c^m r^m{}_{h_1}, \qquad \mu \equiv \mu_2,
    \label{eq:sigma_mu}
\end{equation}
the coordinate choice and Assumption~\ref{assump:unfolding} give
\begin{equation}
    w_c^m r^m{}_{h_2} = 0, \qquad \sigma \neq 0, \qquad \mu \neq 0.
    \label{eq:normalized_unfolding}
\end{equation}
The coefficient $\mu_1$ does not enter the leading-order analysis below. When a model is specified directly in parameters $(h_1,h_2)$ with $\sigma_2 \neq 0$, Eq.~\eqref{eq:param_change} is the map to the coordinates used here. This affects only the labeling of parameter values in numerical examples.

\subsection{Distinguished scaling and reduced dynamics}
\label{sec:scaling}

Truncated at leading order, the vector field on the center manifold is
$\dot A = \sigma h_1 + \mu h_2 A + \nu A^3$. This is quasi-homogeneous: under $A \mapsto \delta A$, $h_2 \mapsto \delta^2 h_2$, $h_1 \mapsto \delta^3 h_1$, every term scales as $\delta^3$. Up to normalization, the weights $(1,2,3)$ of $(A, h_2, h_1)$ are the unique power-law scaling that retains all three terms in a single balance; we fix the weight of $A$ to one. We parametrize this scaling by one small parameter.

\begin{assumption}[Proximity to the organizing center]
\label{assump:scaling}
There is a small parameter $0 < \varepsilon \ll 1$ and fixed $(\tilde{h}_1, \tilde{h}_2) \in \mathbb{R}^2$, independent of $\varepsilon$, such that
\begin{equation}
    h_1 = \varepsilon^{3/2}\, \tilde{h}_1, \qquad h_2 = \varepsilon\, \tilde{h}_2 .
    \label{eq:h_scaling}
\end{equation}
\end{assumption}
With $\delta = \varepsilon^{1/2}$, the state scales as $c \sim \varepsilon^{1/2}$. The parameter $\varepsilon$ measures the weighted distance from the organizing center, and $(\tilde{h}_1,\tilde{h}_2)$ selects a direction in the unfolding. The triple $(\varepsilon, \tilde{h}_1, \tilde{h}_2)$ is redundant under $(\varepsilon, \tilde{h}_1, \tilde{h}_2) \mapsto (\lambda\varepsilon, \lambda^{-3/2}\tilde{h}_1, \lambda^{-1}\tilde{h}_2)$ for $\lambda > 0$. If a unique representative is wanted, it can be fixed by placing $(\tilde{h}_1,\tilde{h}_2)$ on the weighted circle $|\tilde{h}_1|^2 + |\tilde{h}_2|^3 = 1$, i.e.\ by a quasi-homogeneous blow-up of the organizing center. We retain general $(\tilde{h}_1, \tilde{h}_2)$, which allows any fixed direction, including the rays $\tilde{h}_1 = 0$ and $\tilde{h}_2 = 0$.

At this scaling the leading drift is $O(\varepsilon^{3/2})$ in physical variables, and each correction to the truncated normal form is $O(\varepsilon^2)$, a relative factor $\varepsilon^{1/2}$ smaller. The corrections include, for example, the parameter-dependent quadratic term $h_2 A^2$, the cross term $\mu_1 h_1 A$, quartic terms $A^4$ and the second-order constant term $\tfrac12 w_c^m r^m{}_{h_2 h_2} h_2^2$, together with the center-manifold contributions at the same order, such as those from the stable modes displaced by $h_2$.

\begin{assumption}[Controllability of the critical mode]
    Define the control gradient $r^i{}_\alpha = \partial r^i/\partial u^\alpha\big|_{(0,0,0)}$. We assume that $w_c^i r^i{}_\alpha \ne 0$.
    \label{assump:control}
\end{assumption}
This choice ensures that near criticality the control can directly excite the critical mode. We do not consider the broader class of systems in which control enters only through indirect coupling with the stable modes \cite{krener2004control}.

Appending $\dot h = 0$, the extended unforced system has a local center manifold at $(c,h) = (0,0)$, tangent there to $\operatorname{span}\{v_c\} \times \mathbb{R}^2$ \cite{guckenheimer1983nonlinear}. Its slices at fixed small $h$ are one-dimensional, parameter-dependent, locally invariant slow manifolds, on which the stable modes are slaved to the critical amplitude. Tangency to $v_c^i$ is a statement at $h = 0$; for $h \neq 0$ the origin need not be an equilibrium. Under Assumption~\ref{assump:scaling}, the state on this manifold is $c^i \sim \varepsilon^{1/2} A\, v_c^i$. On the slow time $\tau = \varepsilon t$, the amplitude obeys
\begin{equation}
   A' =  \sigma \tilde{h}_1 + \mu \tilde{h}_2\, A + \nu A^3,
   \label{eq:cusp_amplitude}
\end{equation}
as derived in Appendix~\ref{app:weak} with the control set to zero. Equation \eqref{eq:cusp_amplitude} has three equilibria when
\begin{equation}
    \Delta \equiv 27\,\nu^2\sigma^2\tilde{h}_1^2 + 4\,\nu\,\mu^3\tilde{h}_2^3 < 0
    \label{eq:cusp_discriminant}
\end{equation}
and one when $\Delta > 0$; $\Delta$ is minus the discriminant of the cubic $\nu A^3 + \mu\tilde h_2 A + \sigma\tilde h_1$ (equivalently, $-\nu^4$ times the discriminant of the monic cubic obtained by dividing by $\nu$). Within the wedge $\Delta < 0$, the unforced dynamics are bistable for $\nu < 0$.

\begin{assumption}[Locality of initial conditions and target state]\label{assump:locality}
    The initial condition and target state are close to the organizing center:
    \begin{equation}
        |c^i(0)| = O(\varepsilon^{1/2}),
        \qquad
        |c^{*i}| = O(\varepsilon^{1/2}).
    \end{equation}
    Locality of the trajectory itself is assumed separately (Assumption~\ref{assump:existence} and hypothesis (H) of Sec.~\ref{sec:second_order}).
\end{assumption}

\begin{remark}
At the degenerate limit $\varepsilon=0$, Assumption~\ref{assump:scaling} gives $h = 0$ and locality gives $c^i(0)=c^{*i}=0$. The origin is then an equilibrium of the unforced system, so $u^\alpha\equiv 0$ attains $J=0$, with trajectory $c^i\equiv 0$. Since every other control has positive control cost, it is the unique minimizer up to equality almost everywhere.
\end{remark}

\subsection{First-order optimality conditions}
For the objective in Eq.~\eqref{eq:cost_functional}, we introduce the costate variable $p^i(t) \in \mathbb{R}^N$ and write the running Hamiltonian in the usual fashion via Pontryagin's maximum principle \cite{pontryagin2018mathematical,kirk2004optimal}:
\begin{equation}
H[c^i, p^i, u^\alpha] = \frac{1}{2}\left[ \zeta \|c(t) - c^{*}\|^2 + \|u\|^2 \right] + p^i r^i(c^j, u^\alpha, h).
\end{equation}
The first-order necessary optimality conditions for our problem are obtained from the partial derivatives of this function, while the terminal penalty $\frac{\chi}{2}\|c(T)-c^{*}\|^2$ enters through the transversality condition. Differentiation with respect to $p^i$ returns the state dynamics:
\begin{equation}
\frac{d c^i}{d t} = \frac{\partial H}{\partial p^i}  =  r^i(c^j, u^\alpha, h),
\end{equation}
with $c^i(0) = c_0^i$.
The evolution of the costate is governed by
\begin{equation}
\frac{d p^i}{d t} = -\frac{\partial H}{\partial c^i} = -\zeta (c^i - c^{*i}) -  r^j{}_i(c^k, u^\alpha, h)\, p^j,
\end{equation}
subject to the terminal condition
\begin{equation}
p^i(T) = \chi (c^i(T) - c^{*i}).
\end{equation}
The gradient with respect to the control $u^\alpha$ establishes the optimality condition:
\begin{equation}
\frac{\partial H}{\partial u^\alpha} = u^\alpha + r^i{}_{\alpha}(c^j, u^\beta, h)\, p^i = 0,
\end{equation}
which, together with the state and costate equations, constitutes the Hamiltonian boundary-value problem whose solution we analyze asymptotically in Sec.~\ref{sec:asymptotics}. In these conditions $h$ enters only as a fixed parameter. All dependence on $\varepsilon$ arises through Assumptions~\ref{assump:scaling} and \ref{assump:locality} and the scalings of $(\chi, \zeta)$ introduced in Sec.~\ref{sec:asymptotics}.

\begin{assumption}[Existence of a minimizer]\label{assump:existence}
We assume throughout that the problem \eqref{eq:optimization_problem} admits a minimizer $u^\alpha_{\mathrm{opt}} \in \mathcal{U}_{\mathrm{ad}}$ whose trajectory satisfies $\sup_{t\in[0,T]}\|c(t)\| = O(\varepsilon^{1/2})$.
\end{assumption}

\section{Asymptotic analysis of first-order conditions}
\label{sec:asymptotics}
We analyze the optimality conditions perturbatively near the organizing center to obtain leading-order approximations of the extremal solutions. As noted above, at $\varepsilon = 0$ the unique global minimizer is $u^\alpha = 0$ with trajectory $c^i = 0$. We expand about this state. The three regimes of Table~\ref{scal_tab} differ only in the scaling of the tracking penalties, so we treat them together through a single regime index $\omega$ and state the resulting reductions as one theorem (Sec.~\ref{sec:reduction_theorem}). Sec.~\ref{sec:boundary_layers} derives the boundary conditions of the reduced problems, and Sec.~\ref{sec:validation} gives the leading-order solution and feedback law in each case and compares them with numerical optima. The order-by-order derivations for the intermediate and weak regimes are in Appendix~\ref{app:expansions}.

\subsection{Rescaled optimality system}
\label{sec:rescaled_system}

For $\omega \in \{0,1,2\}$, we scale the tracking penalties as
\begin{equation}
    \chi = \varepsilon^{\omega/2}\, \tilde{\chi}, \qquad \zeta = \varepsilon^{\omega}\, \tilde{\zeta},
\label{eq:alpha_rescaling}
\end{equation}
with $\tilde\chi \ge 0$ and $\tilde\zeta > 0$ independent of $\varepsilon$. Since $p^i$ is set by the penalties through the costate equation and terminal condition, we rescale the costate by the same factor as $\chi$:
\begin{equation}
p^i = \varepsilon^{\omega/2}\, q^i .
\label{eq:costate_rescaling}
\end{equation}
Time is rescaled to the scale on which the control acts on the critical mode:
\begin{equation}
\tau = \varepsilon^{\omega/2}\, t, \qquad T_\tau = \varepsilon^{\omega/2}\, T, \qquad \frac{d}{dt} = \varepsilon^{\omega/2} \frac{d}{d\tau}.
\label{eq:time_rescaling}
\end{equation}
The cases $\omega = 0, 1, 2$ are the strong, intermediate and weak regimes of Table~\ref{scal_tab}. Substituting \eqref{eq:alpha_rescaling}--\eqref{eq:time_rescaling} into the optimality conditions of Sec.~\ref{sec:preliminaries}, the state equation changes only through the time derivative, the costate equation is divided by $\varepsilon^{\omega/2}$, and the control condition is written in terms of $q^i$:
\begin{align}
\varepsilon^{\omega/2}\frac{d c^i}{d \tau} &= r^i(c^j, u^\alpha, h), \label{eq:rescaled_state}\\
\varepsilon^{\omega/2}\frac{d q^i}{d \tau} &= -\varepsilon^{\omega/2} \tilde{\zeta}\,(c^i - c^{*i}) - r^j{}_i(c^k, u^\alpha, h)\, q^j, \label{eq:rescaled_costate}\\
 0 &= u^\alpha +    \varepsilon^{\omega/2}\, r^i{}_\alpha(c^j, u^\beta, h)\, q^i, \label{eq:rescaled_control}
\end{align}
with boundary conditions $c^i(0) = c^i_0$ and $q^i(T_\tau) = \tilde\chi\,(c^i(T_\tau) - c^{*i})$. The parameters are $h_1 = \varepsilon^{3/2}\tilde h_1$ and $h_2 = \varepsilon \tilde h_2$ (Assumption~\ref{assump:scaling}), and we write $c^{*i} = \varepsilon^{1/2} c^{*(1),i}$ (Assumption~\ref{assump:locality}).

We expand state, costate and control in half powers of $\varepsilon$:
\begin{align}
c^i &=  \varepsilon^{1/2} c^{(1),i} + \varepsilon c^{(2),i} + \varepsilon^{3/2} c^{(3),i} + O(\varepsilon^2), \label{eq:c_expansion} \\
q^i &=  \varepsilon^{1/2} q^{(1),i} + \varepsilon q^{(2),i} + \varepsilon^{3/2} q^{(3),i} + O(\varepsilon^2), \label{eq:q_expansion} \\
u^\alpha &= \varepsilon^{1/2} u^{(1),\alpha} + \varepsilon u^{(2),\alpha}  + \varepsilon^{3/2} u^{(3),\alpha} + O(\varepsilon^2). \label{eq:u_expansion}
\end{align}
These expansions are understood on bounded sets of the scaled amplitudes, with $\tilde h$, $\tilde\zeta$, $\tilde\chi$ fixed and the slow horizon $T_\tau$ fixed or bounded as $\varepsilon \to 0$. They are not uniform approximations on an $\varepsilon$-independent neighborhood in the unscaled variables $(c, q)$.
Substituting these expansions into $r^i(c^j, u^\alpha, h)$ and collecting powers of $\varepsilon$ gives the expansion \eqref{eq:r_expansion} of Appendix~\ref{app:expansions}. The parameters enter through $h_2 = \varepsilon\tilde h_2$ at $O(\varepsilon)$ and through $h_1 = \varepsilon^{3/2}\tilde h_1$ at $O(\varepsilon^{3/2})$.

The regimes are distinguished by the order at which the control enters the critical amplitude, compared with the drift and with the stable modes. Near the organizing center, the state on the center manifold is $c \sim \varepsilon^{1/2}$, and the intrinsic reduced drift $f(A)$ (Eq.~\eqref{eq:f_def} below) is $O(\varepsilon^{3/2})$ in physical variables. By \eqref{eq:rescaled_control} and \eqref{eq:q_expansion}, the control is $u = O(\varepsilon^{(1+\omega)/2})$. For $\omega = 2$, the control balances the drift, and the reduced problem is fully nonlinear. For $\omega = 1$, the control dominates the drift, and the reduced dynamics are driftless at leading order. For $\omega = 0$, the control is of the same order as the state itself; it then generically excites the stable modes at $O(1)$ relative amplitude on their own $O(1)$ timescale, and no general slow reduction is available.

\subsection{Reduction theorem}
\label{sec:reduction_theorem}

Define the reduced drift
\begin{equation}
    f(A) \equiv \sigma\tilde h_1 + \mu\tilde h_2\, A + \nu A^3,
    \label{eq:f_def}
\end{equation}
with $\sigma$, $\mu$, $\nu$ as in Sec.~\ref{sec:organizing_center}, together with $\beta_c \equiv w_c^i r^i{}_\alpha r^j{}_\alpha w_c^j > 0$ (Assumption~\ref{assump:control}) and the target amplitude $A^* \equiv v_c^i c^{*(1),i}$. The state and costate amplitudes are dual: the state amplitude is the spectral projection $A = w_c^i c^{(1),i}$, while the costate amplitude is $Q = v_c^i q^{(1),i}$, since $q$ is expanded in the left eigenvectors. The target enters only through the costate equation and the terminal condition, so its amplitude is the $v_c$-projection. Since $\|v_c\| = 1$, this is the orthogonal projection of $c^{*(1)}$ onto the critical direction, as dictated by the Euclidean tracking cost. In general $A^* \ne w_c^i c^{*(1),i}$; the two coincide when the target has no stable-mode component, or when the linearization is symmetric, so that $w_c = v_c$.

\begin{theorem}[Reduction of the optimality system]
\label{thm:reduction}
Let Assumptions~\ref{assump:bifurcation}--\ref{assump:locality} hold, and scale the penalties as in \eqref{eq:alpha_rescaling} with $\omega \in \{0,1,2\}$.

(i) In the strong regime, $\omega = 0$, at $\varepsilon = 0$ the linearization of the optimality system about the organizing center has no eigenvalues on the imaginary axis. The leading-order problem, at $O(\varepsilon^{1/2})$, is the linear--quadratic tracking problem for the linearized system \eqref{eq:strong_state}--\eqref{eq:strong_control} on the full state space. There is no center subspace and no $\varepsilon$-induced separation of timescales, so no general cusp-based slow reduction is available.

(ii) In the intermediate regime, $\omega = 1$, for sufficiently small $\varepsilon > 0$, the rescaled optimality system has a two-dimensional locally invariant manifold $\mathcal{M}_\varepsilon$ in a neighborhood of the organizing center. At $\varepsilon = 0$ it is tangent to $\operatorname{span}\{(v_c, 0), (0, w_c)\}$ in $(c,q)$-space. On $\mathcal{M}_\varepsilon$,
\begin{equation}
c^i = \varepsilon^{1/2}\big(A\, v_c^i + O(\varepsilon^{1/2})\big), \qquad
q^i = \varepsilon^{1/2}\big(Q\, w_c^i + O(\varepsilon^{1/2})\big),
\label{eq:manifold_param}
\end{equation}
and the amplitudes obey
\begin{equation}
    A' = -\beta_c Q + O(\varepsilon^{1/2}), \qquad
    Q' = -\tilde\zeta\,(A - A^*) + O(\varepsilon^{1/2}).
    \label{eq:thm_intermediate}
\end{equation}
(iii) In the weak regime, $\omega = 2$, the conclusions of (ii) hold with the reduced dynamics
\begin{equation}
    A' = f(A) - \beta_c Q + O(\varepsilon^{1/2}), \qquad
    Q' = -\tilde\zeta\,(A - A^*) - f'(A)\, Q + O(\varepsilon^{1/2}).
    \label{eq:thm_weak}
\end{equation}

In cases (ii) and (iii), the manifold $\mathcal{M}_\varepsilon$ is not attracting. In the limit $\varepsilon \to 0$, its transverse spectrum at the organizing center, on the fast time $t$, consists of the $N-1$ stable eigenvalues $\lambda_s$ of the state and the $N-1$ unstable eigenvalues $-\lambda_s$ of the costate. For $\varepsilon > 0$ these are perturbed, and on the slow time $\tau$ the transverse rates are multiplied by $\varepsilon^{-\omega/2}$. The reduced flow is Hamiltonian in $(A,Q)$ with
\begin{equation}
    H_\omega(A,Q) = \frac{\tilde\zeta}{2}(A - A^*)^2 - \frac{\beta_c}{2}Q^2 + \varepsilon^{(2-\omega)/2}\, Q\, f(A),
    \label{eq:thm_hamiltonian}
\end{equation}
up to the errors stated in \eqref{eq:thm_intermediate} and \eqref{eq:thm_weak}.
\end{theorem}

\begin{proof}
Throughout we work in the fast time $t$ and write $\delta = \varepsilon^{1/2}$. In this notation, Eqs.~\eqref{eq:rescaled_state}--\eqref{eq:rescaled_control} read
\begin{equation*}
\frac{dc^i}{dt} = r^i(c, u, h), \qquad
\frac{dq^i}{dt} = -\delta^{\omega}\tilde\zeta\,(c^i - \delta\, c^{*(1),i}) - r^j{}_i(c,u,h)\, q^j, \qquad
u^\alpha = -\delta^{\omega}\, r^i{}_\alpha(c,u,h)\, q^i,
\end{equation*}
with $h_1 = \delta^3\tilde h_1$ and $h_2 = \delta^2\tilde h_2$. By the implicit function theorem, the control condition has a unique small solution $u = u(c, q, \delta, h)$ near the origin, smooth in its arguments, so the system is autonomous and smooth in $(c, q, \delta)$.

Consider first case (i). For $\omega = 0$, at $\delta = 0$, the linearization about $(c,q) = (0,0)$ is the Hamiltonian matrix
\begin{equation*}
\mathcal{H} = \begin{pmatrix} R & -BB^\top \\ -\tilde\zeta\, I & -R^\top \end{pmatrix}, \qquad R_{ij} = r^i{}_j, \quad B_{i\alpha} = r^i{}_\alpha .
\end{equation*}
The only eigenvalue of $R$ with nonnegative real part is the simple eigenvalue $\lambda_c = 0$, with left eigenvector $w_c$, and $w_c^\top B \neq 0$ by Assumption~\ref{assump:control}. By the Popov--Belevitch--Hautus test, $(R, B)$ is therefore stabilizable, and since $\tilde\zeta > 0$, $(R, \tilde\zeta^{1/2} I)$ is observable. Hence $\mathcal{H}$ has no eigenvalues on the imaginary axis, and its spectrum splits into $N$ eigenvalues with negative and $N$ with positive real part \cite{abou2012matrix}. The zero eigenvalue of $R$ therefore does not persist: through the $O(1)$ coupling $-BB^\top$, the critical direction relaxes on the same $O(1)$ timescale as the stable modes, and there is no slow amplitude. Since $(c,q) = O(\delta)$ and the nonlinear terms are $O(\delta^2)$, the leading-order problem at $O(\delta)$ is the linear--quadratic problem of Sec.~\ref{sec:strong_regime}, which is well posed because it is strictly convex for $\tilde\zeta > 0$, $\tilde\chi \ge 0$.

For cases (ii) and (iii), append the trivial equation $d\delta/dt = 0$. Every term that couples $c$ and $q$, and every term containing $h$, $c^*$ or $u$, carries at least one factor of $\delta$. All such terms are therefore nonlinear in the extended variables $(c, q, \delta)$. The linearization about the origin is consequently block diagonal, $\operatorname{diag}(R, -R^\top, 0)$. Its spectrum consists of a zero eigenvalue of multiplicity three, with eigenvectors $(v_c, 0, 0)$, $(0, w_c, 0)$ and $e_{\delta}$; the $N-1$ eigenvalues $\lambda_s$ of $R$, with $\mathrm{Re}\,\lambda_s < 0$; and the $N-1$ eigenvalues $-\lambda_s$ of $-R^\top$, with positive real part.

The hyperbolic part contains both stable and unstable directions. The center manifold theorem in its general form, which does not require the hyperbolic part to be stable, gives a local $C^k$ center manifold $W^c$ of dimension three \cite{kelley1967stable, guckenheimer1983nonlinear, vanderbauwhede1989centre}. It is locally invariant and tangent to the center subspace. It contains every solution that remains in a sufficiently small neighborhood of the origin for all $t \in \mathbb{R}$ \cite{vanderbauwhede1989centre, haragus2011local}. Since $\delta$ is a constant of motion, $W^c$ is foliated by the two-dimensional invariant slices $\mathcal{M}_\varepsilon = W^c \cap \{\delta = \varepsilon^{1/2}\}$. The manifold $W^c$ is not unique, but its Taylor expansion at the origin is \cite{carr1981applications}. The reduced equations below are therefore independent of the choice of $W^c$ at any finite order. Because the transverse directions include the unstable costate modes, $W^c$ is not attracting, and the reduction principle of the dissipative case does not apply.

To compute the reduced vector field on $\mathcal{M}_\varepsilon$, write $c = \delta(A v_c + \phi)$ and $q = \delta(Q w_c + \psi)$. Here $\phi$ lies in the stable subspace $\operatorname{span}\{v_s\}$ and $\psi$ in the corresponding adjoint subspace $\operatorname{span}\{\overline{w_s}\}$, and both are $O(\delta)$. Projecting the state equation onto $w_c$ and the costate equation onto $v_c$ gives
\begin{align*}
\frac{dA}{dt} &= \delta^2 f(A) - \delta^{\omega}\, \beta_c Q + O(\delta^3 + \delta^{1+\omega}), \\
\frac{dQ}{dt} &= -\delta^{\omega} \tilde\zeta\,(A - A^*) - \delta^2 f'(A)\, Q + O(\delta^3 + \delta^{1+\omega}).
\end{align*}
These projections are computed order by order in Appendices~\ref{app:intermediate} and \ref{app:weak}. Three terms of lower order are absent. The $O(\delta)$ terms proportional to $w_c^i r^i{}_{jk} v_c^j v_c^k$, namely $A^2$ in the first equation and $AQ$ in the second, vanish by Assumption~\ref{assump:cusp}. The $O(\delta)$ constant $\delta\, w_c^i r^i{}_{h_2}\tilde h_2$ in $dA/dt$, arising from $h_2 = \delta^2\tilde h_2$, vanishes in the parameter coordinates \eqref{eq:normalized_unfolding}.
Changing to $\tau = \delta^{\omega} t$ and dividing by $\delta^{\omega}$ multiplies the drift terms by $\delta^{2-\omega} = \varepsilon^{(2-\omega)/2}$. For $\omega = 1$ this factor is $O(\varepsilon^{1/2})$, which gives \eqref{eq:thm_intermediate}; for $\omega = 2$ it equals one, which gives \eqref{eq:thm_weak}.

Finally, after elimination of $u$, the Pontryagin system is Hamiltonian in $(c, p)$, with Hamiltonian $\bar H(c,p) = H(c,p,u(c,p))$, the Hamiltonian evaluated on this local stationary branch of $H_u = 0$, and symplectic form $dc^i \wedge dp^i$. The center manifold of a Hamiltonian system is a symplectic submanifold, and the reduced flow is Hamiltonian with respect to the restricted form and the restricted Hamiltonian \cite{mielke1991hamiltonian}. On $\mathcal{M}_\varepsilon$, using $p = \varepsilon^{\omega/2} q$ and \eqref{eq:manifold_param},
\begin{equation*}
dc^i\wedge dp^i\big|_{\mathcal{M}_\varepsilon} = \varepsilon^{(2+\omega)/2}\,\big(1 + O(\varepsilon^{1/2})\big)\, dA \wedge dQ,
\qquad
\bar H\big|_{\mathcal{M}_\varepsilon} = \varepsilon^{1+\omega}\, H_\omega(A,Q) + \text{const} + \text{h.o.t.}
\end{equation*}
In the second relation, the additive constant comes from the component of the target orthogonal to $v_c$. The term $Q f(A)$ arises from $p^i r^i$, and the $Q^2$ term from $-\tfrac12 \|r^i{}_\alpha p^i\|^2$. Hamilton's equations for $H_\omega$ in the slow time $\tau$ reproduce \eqref{eq:thm_intermediate} and \eqref{eq:thm_weak}.
\end{proof}

\begin{remark}[Rigorous, formal and numerical results]
\label{rem:optimal_trajectories}
The results of this section rest on three levels of justification. Theorem~\ref{thm:reduction} is rigorous. It establishes the existence of $\mathcal{M}_\varepsilon$, its saddle-type transverse spectrum, and the Hamiltonian reduced flow on it. It is a statement about the invariant manifold, not about solutions of the boundary-value problem. By the center-manifold property, a solution that remains in a small neighborhood of the organizing center for all $t \in \mathbb{R}$ lies on $\mathcal{M}_\varepsilon$. A finite-horizon extremal with generic boundary data is not of this type: its initial state has components along the stable directions, and its terminal costate has components along the unstable ones. Sec.~\ref{sec:boundary_layers} derives the structure of such an extremal by formal matched asymptotics. In that construction, the extremal consists of an initial layer, an arc near $\mathcal{M}_\varepsilon$ governed at leading order by the reduced flow, and a terminal layer, and the reduced problem inherits the projected boundary conditions \eqref{eq:reduced_bcs}. We do not prove error bounds for this construction, so the leading-order solutions of Sec.~\ref{sec:validation} are asymptotic approximations of extremals. Their agreement with numerical optima of the full problem, including the $O(\varepsilon)$ decay of the trajectory error, is observed numerically for the two testbeds of Appendix~\ref{app:numerics}. Whether these extremals are local minimizers is a separate question. Sec.~\ref{sec:second_order} addresses it under a hypothesis that encodes the formal layer structure.
\end{remark}

\begin{remark}[Dependence on the unfolding]
In cases (i) and (ii), the unfolding parameters do not enter the leading-order problem. In case (iii), they enter only through $f(A)$, as the constant $\sigma\tilde h_1$ and the linear coefficient $\mu\tilde h_2$. The parameter coordinates of Sec.~\ref{sec:organizing_center} matter here. In other coordinates, a constant $\sigma_2\tilde h_2$ would enter the critical amplitude at $O(\varepsilon)$, which is comparable to the control for $\omega = 1$ and dominant for $\omega = 2$, so both reduced problems would be altered.
\end{remark}

\subsection{Boundary conditions and boundary layers}
\label{sec:boundary_layers}

The optimality system carries $2N$ scalar boundary conditions, $c^i(0) = c^i_0$ and
$p^i(T) = \chi(c^i(T) - c^{*i})$. In the strong regime ($\omega = 0$) no slow reduction is available in general, and the
leading-order problem is the full $2N$-dimensional system. For $\omega \in \{1,2\}$ the leading-order flow is the
two-dimensional $(A,Q)$ system on $\mathcal{M}_\varepsilon$ and can accept only two boundary conditions. The remaining $2(N-1)$
conditions are absorbed in boundary layers, which we identify here. As anticipated in Remark~\ref{rem:optimal_trajectories}, these layers are attached to the stable directions of $\mathcal{M}_\varepsilon$ at the initial time and to its unstable directions at the terminal time.

Write $\tau = \varepsilon^{\omega/2} t$, and let $t$ and $\hat{t} = T - t$ be the fast variables at
the two ends. Decompose the solution into outer and layer parts,
$c^i = c^i_{\mathrm{out}}(\tau) + \hat{c}^i(t)$ and
$p^i = p^i_{\mathrm{out}}(\tau) + \hat{p}^i(\hat{t})$. For $\omega \in \{1,2\}$, the control is
smaller than the state by a factor $\varepsilon^{\omega/2}$, and the tracking term
$\zeta(c^i - c^{*i})$ is smaller than the costate by the same factor. At leading order,
the layer equations therefore decouple into a linear pair,
\begin{equation}
    \frac{d \hat{c}^i}{d t} = r^i{}_j\, \hat{c}^j,
    \qquad
    \frac{d \hat{p}^i}{d \hat{t}} = r^j{}_i\, \hat{p}^j ,
\label{eq:layer_equations}
\end{equation}
with relative corrections of $O(\varepsilon^{1/2})$ in both regimes, since the cusp condition removes only the critical projection of the quadratic state terms. By Assumption~\ref{assump:bifurcation}, the spectrum of $r^i{}_j$ off the critical direction satisfies
$\mathrm{Re}\,\lambda_s \le -\Lambda$ with $\Lambda = \min_s |\mathrm{Re}\,\lambda_s| > 0$, so
\begin{equation}
    \hat{c}^i(t) = \sum_s \hat{a}_s\, e^{\lambda_s t}\, v_s^i ,
    \qquad
    \hat{p}^i(\hat{t}) = \sum_s \hat{b}_s\, e^{\lambda_s \hat{t}}\, \overline{w_s}^{\,i} .
\label{eq:layer_solutions}
\end{equation}
Both families decay at least at rate $\Lambda$ on the fast scale, that is over
$\tau = O(\varepsilon^{\omega/2})$, and are negligible beyond
$\tau \sim \Lambda^{-1}\varepsilon^{\omega/2}\ln(1/\varepsilon)$. The state relaxes onto the
manifold in an initial layer, and the costate relaxes onto it, backward in time, in
a terminal layer. A terminal state layer or an initial costate layer would grow into the
interior and is inadmissible. The two layers are transverse to the critical directions: by biorthogonality,
$w_c^i \hat{c}^i = 0$ and $v_c^i \hat{p}^i = 0$. The critical projections of the boundary
data are therefore unaffected at leading order. Projecting the initial condition onto $w_c$
and the terminal costate condition onto $v_c$ gives the reduced boundary conditions
\begin{equation}
    A(0) = \varepsilon^{-1/2} w_c^i c^i_0 ,
    \qquad
    Q(T_\tau) = \tilde{\chi}\big(A(T_\tau) - A^*\big),
\label{eq:reduced_bcs}
\end{equation}
for $\omega \in \{1,2\}$. The layers re-enter the slow dynamics only through higher-order corrections. These include the quadratic coupling
$w_c^i r^i{}_{jk} \hat{c}^j \hat{c}^k$ and its costate analog. There are also linear couplings. The terminal costate layer forces the critical state through the control, $-w_c^\top \mat{B}\mat{B}^\top \hat{\vec p}$, which need not vanish when $v_c^i \hat p^i = 0$. For nonorthogonal eigenvectors, the initial state layer enters the critical costate through the tracking term. Since the layers act over a time
$t = O(1)$ and $\hat{c}^i = O(\varepsilon^{1/2})$ for data admitted by
Assumption~\ref{assump:locality}, these couplings displace the reduced boundary data \eqref{eq:reduced_bcs} by a
relative $O(\varepsilon^{1/2})$, which is the order of the first term
neglected in \eqref{eq:c_expansion}. Boundary data consistent beyond leading order
would require the conditions to be
projected along the isochrons of the manifold rather than along the spectral
projection, in the manner of Ref.~\cite{roberts1989appropriate}. We use the
leading-order projection \eqref{eq:reduced_bcs} throughout and do not
pursue this correction.

\subsection{Leading-order solutions and numerical validation}
\label{sec:validation}

We now give the leading-order solution in each regime, in feedback form where possible, and compare it with numerical optima of the full problem. The comparisons use two testbeds detailed in Appendix~\ref{app:numerics}: a bistable-switch kinetic model and a Turing-pattern reaction--diffusion system. For $\omega \in \{1,2\}$, the reduced problems are posed with the boundary conditions \eqref{eq:reduced_bcs}.

\subsubsection{Case (i): strong control, $\omega = 0$}
\label{sec:strong_regime}

With $\omega = 0$, the rescalings \eqref{eq:alpha_rescaling}--\eqref{eq:time_rescaling} are trivial, and $q^i = p^i$, $\tau = t$. The lowest nontrivial order is $O(\varepsilon^{1/2})$:
\begin{align}
    \frac{d c^{(1),i}}{dt } &= r^i{}_j c^{(1), j} + r^i{}_\alpha u^{(1), \alpha}, \label{eq:strong_state}\\
    -\frac{d p^{(1), i}}{dt } &= \zeta(c^{(1), i} - c^{*(1), i}) + r^j{}_i p^{(1),j}, \label{eq:strong_costate}\\
    u^{(1), \alpha} &= - r^i{}_\alpha p^{(1),i}. \label{eq:strong_control}
\end{align}
These are the optimality conditions of an LQR tracking problem for the linear system with state matrix $r^i{}_j$. The parameters $h$ first enter at $O(\varepsilon)$ and so do not appear. As in any linear--quadratic tracking problem, the costate is affine in the state, $\vec{p}^{(1)} = \mat{K}(t)\vec{c}^{(1)} + \vec{b}(t)$, where $\mat{K}\colon [0, T] \to M_{N\times N}(\mathbb{R})$ solves the matrix Riccati equation
\begin{equation}
    - \frac{d \mat{K}}{dt} = \mat{R}^\top \mat{K} +   \mat{K} \mat{R} - \mat{K} \mat{B} \mat{B}^\top \mat{K} + \zeta \mat{I} \label{eq:riccati_strong}
\end{equation}
with $R_{ij} = r^i{}_j$, $B_{i\alpha} = r^i{}_\alpha$ and boundary condition $\mat{K}(T) = \chi \mat{I}$. The offset $\vec b$ solves
\begin{equation}
    -\frac{d \vec{b}}{dt} = (\mat{R}^\top - \mat{K} \mat{B} \mat{B}^\top)\vec{b} - \zeta \vec{c}^{*(1)}
\end{equation}
with $\vec{b}(T) = - \chi \vec{c}^{*(1)}$. The gain can be computed numerically, which gives a practical feedback law.

\begin{figure}[!t]
    \centering
    \includegraphics{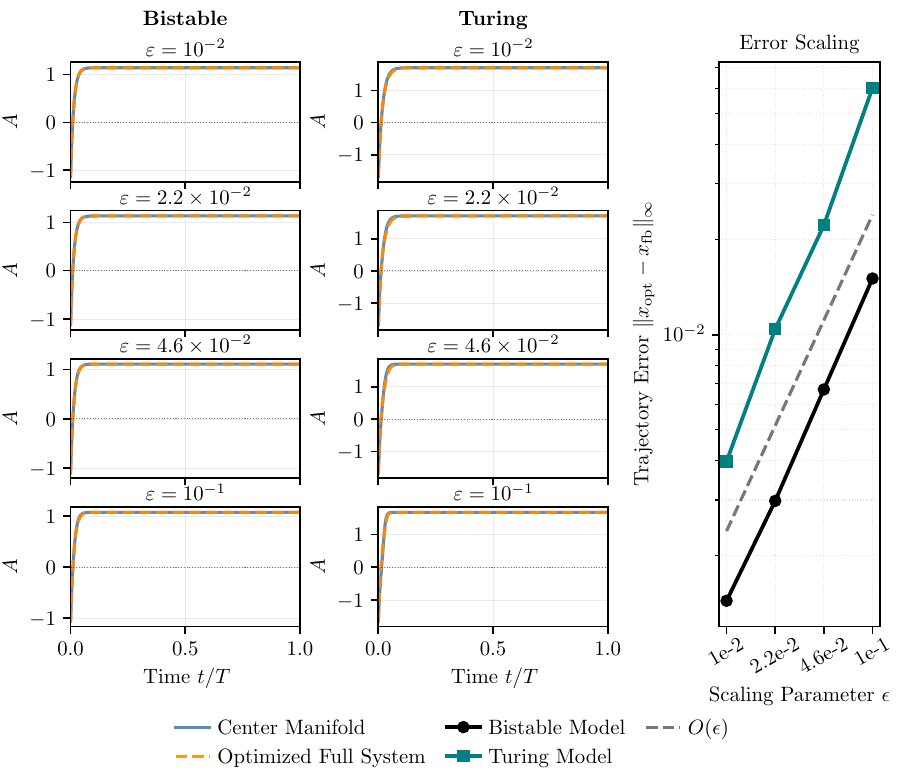}
    \caption{Strong control regime validation, combining the bistable-switch (left column) and Turing-pattern (center column) testbeds. Each row shows the leading critical-mode amplitude $A(t)$ at a fixed $\varepsilon$ (increasing top to bottom over $\varepsilon\in[10^{-2},10^{-1}]$), comparing the asymptotic LQR feedback law against the numerically optimized full system. Right column: the full-state $L_\infty$ trajectory error between the two solutions as a function of $\varepsilon$ for both testbeds, together with a reference $O(\varepsilon)$ slope.}
    \label{fig:lqr}
\end{figure}

Figure~\ref{fig:lqr} compares the asymptotic LQR solution with the numerical optimum of the full nonlinear problem; the error decreases as $O(\varepsilon)$.

This regime is where the universal center-manifold structure is lost. Consistent with Theorem~\ref{thm:reduction}(i), the term $\mat{K}\mat{B}\mat{B}^\top\mat{K}$ in \eqref{eq:riccati_strong} generically couples the critical and stable modes at leading order: unless the control acts only along a critical direction decoupled from the stable modes, it excites the stable modes to finite amplitude despite their intrinsic decay.

\subsubsection{Case (ii): intermediate control, $\omega = 1$}
\label{sec:intermediate_regime}

With $\omega = 1$, we have $\chi = \varepsilon^{1/2}\tilde\chi$, $\zeta = \varepsilon\tilde\zeta$, $p^i = \varepsilon^{1/2}q^i$ and $\tau = \varepsilon^{1/2}t$. The leading-order solution is $c^{(1),i} = A(\tau) v_c^i$, $q^{(1),i} = Q(\tau) w_c^i$, and solvability at $O(\varepsilon)$ fixes the amplitudes (Appendix~\ref{app:intermediate}):
\begin{equation}
    A' = - \beta_{c} Q, \qquad
    Q' = - \tilde{\zeta}(A  - A^*),
    \label{eq:2nd_opt}
\end{equation}
which is case (ii) of Theorem~\ref{thm:reduction}. The unfolding first enters at $O(\varepsilon^{3/2})$, as a correction of relative size $O(\varepsilon^{1/2})$.

\begin{figure}[!t]
    \centering
    \includegraphics{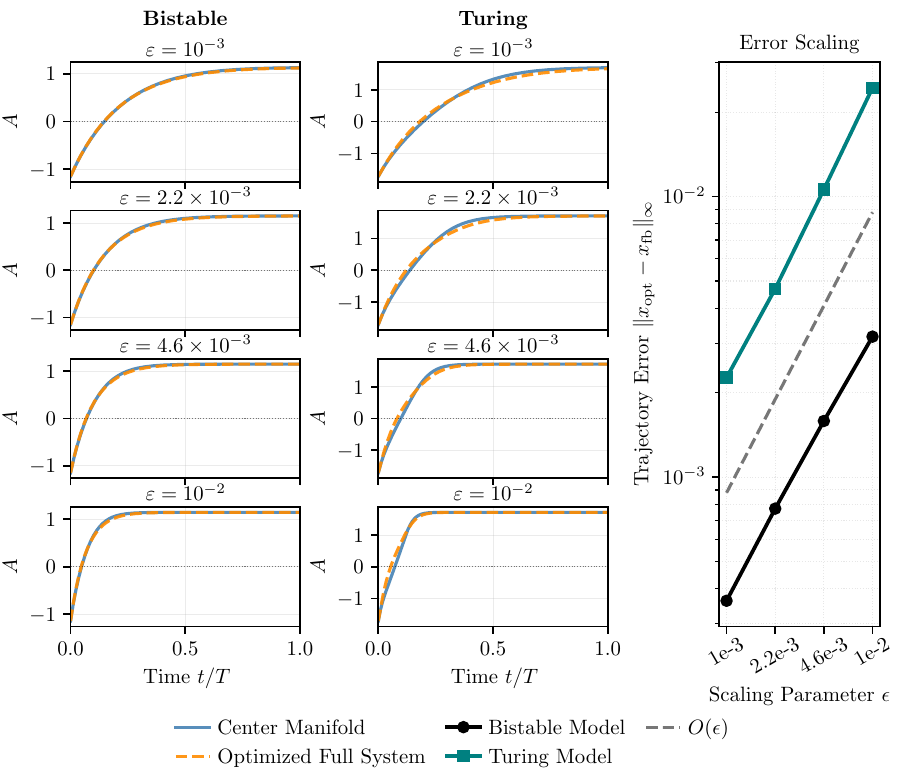}
    \caption{Intermediate control regime validation, combining the bistable-switch (left column) and Turing-pattern (center column) testbeds. Each row shows the leading critical-mode amplitude $A(t)$ at a fixed $\varepsilon$ (increasing top to bottom over $\varepsilon\in[10^{-3},10^{-2}]$), comparing the tanh-Riccati asymptotic feedback law against the numerically optimized full system. Right column: the full-state $L_\infty$ trajectory error between the two solutions as a function of $\varepsilon$ for both testbeds, together with a reference $O(\varepsilon)$ slope.}
    \label{fig:feedback}
\end{figure}

These equations describe the optimal control of a one-dimensional, driftless system. For constant $A^*$, they admit the linear feedback form $Q = \Pi(\tau) (A - A^*)$, where $\Pi$ solves the scalar Riccati equation
\begin{equation}
    - \frac{d \Pi}{d \tau} = \tilde{\zeta} - \beta_c \Pi^2,  \label{eq:scalar_riccati}
\end{equation}
with boundary condition $\Pi(T_\tau) = \tilde{\chi}$. With $k = \sqrt{\tilde\zeta\beta_c}$ and $\hat\chi = \tilde\chi\sqrt{\beta_c/\tilde\zeta}$, the gain is
\begin{equation}
    \Pi(\tau) = \sqrt{\frac{\tilde{\zeta}}{\beta_c}}\;
    \frac{\hat\chi + \tanh\!\big(k(T_\tau - \tau)\big)}{1 + \hat\chi\tanh\!\big(k(T_\tau - \tau)\big)},
\end{equation}
and the leading-order feedback law is
\begin{equation}
    u^{(2), \alpha}(\tau) = - \Pi(\tau)\, r^i{}_\alpha w_c^i\,\big(A - A^*\big), \qquad A = w_c^j c^{(1),j}, \quad A^* = v_c^j c^{*(1),j}.
\end{equation}
The two projections differ: the state amplitude is read off spectrally, with $w_c$, and the target amplitude orthogonally, with $v_c$.
Figure~\ref{fig:feedback} compares this feedback law with numerical optima; it captures the dominant control and state behavior near criticality. The stable-mode corrections and the construction of higher orders are given in Appendix~\ref{app:intermediate}.

\subsubsection{Case (iii): weak control, $\omega = 2$}
\label{sec:weak_regime}

With $\omega = 2$, we have $\chi = \varepsilon\tilde\chi$, $\zeta = \varepsilon^2\tilde\zeta$, $p^i = \varepsilon q^i$ and $\tau = \varepsilon t$. The leading-order solution again has the form $c^{(1),i} = A(\tau) v_c^i$, $q^{(1),i} = Q(\tau) w_c^i$, but the control is now absent through $O(\varepsilon)$, and solvability enters only at $O(\varepsilon^{3/2})$ (Appendix~\ref{app:weak}). The amplitudes obey
\begin{align}
    A'(\tau) &= f(A) -  \beta_{c} Q = \sigma\tilde h_1 + \mu\tilde h_2\, A  + \nu A^3 -  \beta_{c} Q,  \label{eq:Aweak}\\
    Q'(\tau) &= \tilde{\zeta}(A^* - A(\tau)) - f'(A)\, Q = \tilde{\zeta}(A^* - A(\tau)) - \big(\mu\tilde h_2 + 3 \nu A^2\big) Q, \label{eq:Qweak}
\end{align}
which is case (iii) of Theorem~\ref{thm:reduction}. The reduced dynamics are now fully nonlinear. Equations \eqref{eq:Aweak}--\eqref{eq:Qweak} are the first-order necessary conditions for minimizing
\begin{equation}
    \mathcal{L} = \frac{\tilde\chi}{2}\big(A(T_\tau) - A^*\big)^2 + \frac{1}{2} \int_0^{T_\tau} \big(\tilde{\zeta}(A - A^*)^2 + v^2\big)\, d\tau
\end{equation}
subject to the state dynamics $A'(\tau) =  f(A) + \sqrt{\beta_{c}}\, v$, with costate $Q$ and $v = -\sqrt{\beta_c}\, Q$.

\begin{figure}[!t]
    \centering
    \includegraphics{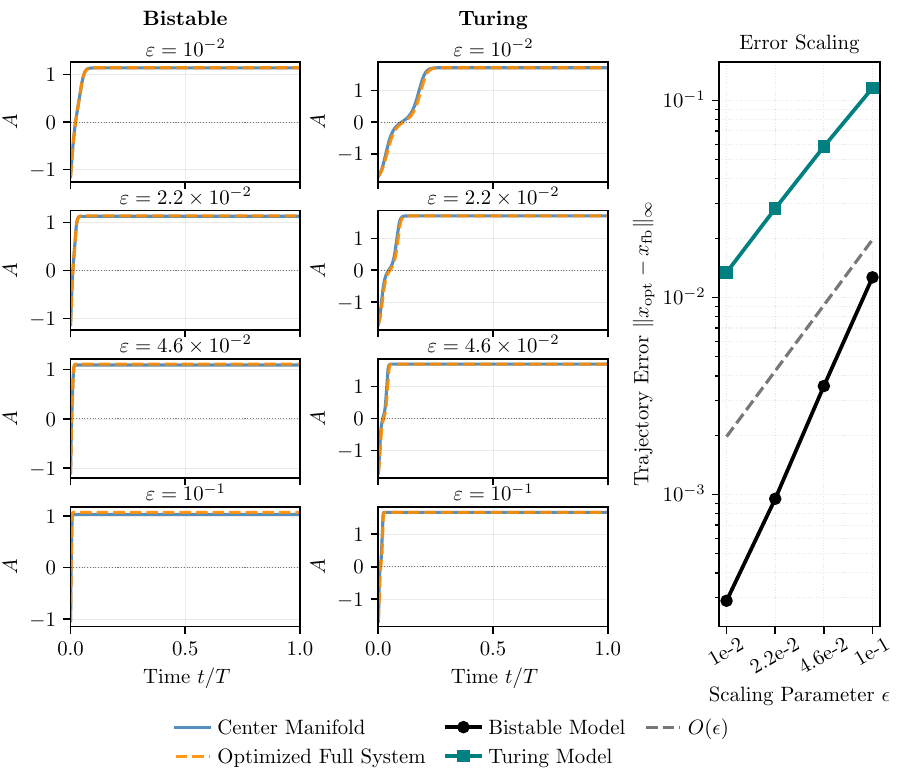}
    \caption{Weak control regime validation, combining the bistable-switch (left column) and Turing-pattern (center column) testbeds. Each row shows the leading critical-mode amplitude $A(t)$ at a fixed $\varepsilon$ (increasing top to bottom over $\varepsilon\in[10^{-2},10^{-1}]$), comparing the asymptotic solution of the reduced problem against the numerically optimized full system. Right column: the full-state $L_\infty$ trajectory error between the two solutions as a function of $\varepsilon$ for both testbeds, together with a reference $O(\varepsilon)$ slope.}
    \label{fig:nonlinear}
\end{figure}

The solution can also be written in feedback form using the reduced Hamiltonian, $H_2$ of \eqref{eq:thm_hamiltonian}:
\begin{equation}
    H_{\text{reduced}} = \frac{\tilde{\zeta}}{2} (A - A^*)^2 - \frac{\beta_c Q^2 }{2} + Q f(A).
\end{equation}
For constant $A^*$, $H_{\text{reduced}}$ does not depend explicitly on time and is conserved along extremals. Let $E$ be its value on a selected extremal with initial amplitude $A_0 = A(0)$ and horizon $T_\tau$; since the problem is nonconvex, several extremals with the same $A_0$ and $T_\tau$ may have different energies. Since $H_{\text{reduced}}$ is quadratic in $Q$, the level set $H_{\text{reduced}} = E$ gives two feedback branches,
\begin{equation}
   Q_\pm(A, E) =  \frac{f(A) \pm \sqrt{f(A)^2+ \beta_c  \left(\tilde{\zeta}(A-A^*)^2 -2E(A_0, T_\tau )\right) } }{\beta_c }.
\end{equation}
Which branch applies depends on the initial data and the target, and the extremal switches branch at each turning point of $A$. The feedback law therefore need not be single-valued in $A$, for example when the trajectory overshoots. These formulas describe only the selected extremal; they do not define a global feedback law. The level $E(A_0, T_\tau)$ is fixed implicitly by the boundary conditions, and we do not have a closed-form expression for it. Since $p^i=\varepsilon q^i$ and $q^{(1),i}=Qw_c^i$ enters at $O(\varepsilon^{1/2})$, we have $p^i=\varepsilon^{3/2}Qw_c^i+O(\varepsilon^2)$, and the physical control is
\begin{align}
    u^{\alpha} &= \varepsilon^{3/2} u^{(3), \alpha} + O(\varepsilon^2) = -\varepsilon^{3/2} Q\, r^i{}_\alpha w_c^i + O(\varepsilon^2) \nonumber \\
    &= - \frac{\varepsilon^{3/2} }{\beta_c}  \left[f(A) \pm\sqrt{f(A)^2 + \beta_c \left(\tilde{\zeta}(A-A^*)^2 - 2E(A_0, T_\tau)\right)} \right] r^i{}_\alpha w_c^i + O(\varepsilon^2).
\end{align}
Figure~\ref{fig:nonlinear} compares the weak-regime reduction with numerical optima. For the bistable switch, the full system is driven by this level-set feedback law; for the Turing model, by the solution of the reduced boundary-value problem (Appendix~\ref{app:numerics}).

\section{Second-order conditions}
\label{sec:second_order}

The extremals of Sec.~\ref{sec:asymptotics} satisfy the first-order necessary conditions only. In this section we derive the leading-order second variation in each regime of Theorem~\ref{thm:reduction}, and state when it is coercive. For the weak regime we then give focal-point tests that use only the extremal itself, and apply them to its layered extremals in Sec.~\ref{sec:boundary_layer}.

Let $(\bar c^i, \bar p^i, \bar u^\alpha)$ solve the first-order optimality conditions of Sec.~\ref{sec:preliminaries}. For the cost \eqref{eq:cost_functional}, the second variation with respect to a control perturbation $\tilde u^\alpha$ is
\begin{equation}
\delta^2 J[\tilde{u}]
= \chi\,\|\tilde{c}(T)\|^2
  + \int_0^T \Big[
      \zeta\,\|\tilde{c}\|^2
    + \|\tilde{u}\|^2
    + \bar p^m \big(
        r^m{}_{ij}\tilde{c}^i\tilde{c}^j
        + 2r^m{}_{i\alpha}\tilde{c}^i\tilde{u}^\alpha
        + r^m{}_{\alpha\beta}\tilde{u}^\alpha\tilde{u}^\beta
      \big)
    \Big] dt,
\label{eq:second_variation}
\end{equation}
so that $J[\bar u+\tilde u]=J[\bar u]+\tfrac12\delta^2J[\tilde u]+o(\|\tilde u\|^2_{L^2})$ as $\|\tilde u\|_{L^\infty}\to0$, where $\tilde c^i$ solves the linearization about the extremal,
\begin{equation}
    \frac{d\tilde{c}^i}{dt} = r^i{}_j(\bar c, \bar u, h)\, \tilde{c}^j + r^i{}_\alpha(\bar c, \bar u, h)\, \tilde{u}^\alpha, \qquad \tilde{c}^i(0) = 0,
\label{eq:linearized_state}
\end{equation}
and all derivatives of $r^m$ are evaluated along $(\bar c, \bar u, h)$. When $r$ is nonlinear in $u$, the expansion holds only for $L^\infty$-small perturbations, although the remainder is measured in $L^2$ (the two-norm discrepancy \cite{maurer1981first}). Local optimality below is therefore understood in the weak sense: minimality in an $L^\infty$ neighborhood of $\bar u$, with growth measured in $L^2$. The Legendre--Clebsch matrix is $H_{u^\alpha u^\beta} = \delta_{\alpha\beta} + \bar p^m r^m{}_{\alpha\beta}$. Since $\bar p = O(\varepsilon^{1/2})$ in every regime, it is positive definite for small $\varepsilon$, and the strengthened Legendre condition holds throughout.

In the intermediate and weak regimes, the second variation reduces to a scalar quadratic form, which we now define. For $\omega \in \{1,2\}$ let $\theta_\omega = \lim_{\varepsilon\to0}\varepsilon^{(2-\omega)/2}$, so $\theta_1 = 0$ and $\theta_2 = 1$. Let $(\bar A, \bar Q)$ be the reduced extremal of Theorem~\ref{thm:reduction}, with $f$ the cusp drift \eqref{eq:f_def}. Define
\begin{equation}
    \Theta(\tau) = \theta_\omega\, f'(\bar A), \qquad
    H_{AA}(\tau) = \tilde\zeta + \theta_\omega\, f''(\bar A)\,\bar Q = \tilde\zeta + 6\theta_\omega\nu\bar A\bar Q .
    \label{eq:ell_HAA}
\end{equation}
$H_{AA}$ is $\partial^2 H_\omega/\partial A^2$ for the reduced Hamiltonian \eqref{eq:thm_hamiltonian}, and $\Theta$ is the coefficient of the linearized reduced state equation. For a scalar perturbation $v \in L^2(0,T_\tau)$ let $\xi$ solve
\begin{equation}
    \xi' = \Theta\,\xi + \sqrt{\beta_c}\,v, \qquad \xi(0) = 0,
    \label{eq:xi_lin}
\end{equation}
and set
\begin{equation}
    \mathcal{Q}_\omega[v] = \tilde\chi\,\xi(T_\tau)^2 + \int_0^{T_\tau}\big[v^2 + H_{AA}\,\xi^2\big]\,d\tau .
    \label{eq:2var_weak}
\end{equation}
This is the second variation of the reduced problem whose extremals are given by Theorem~\ref{thm:reduction}(ii)--(iii). For $\omega=2$ it is the second variation of the cost $\mathcal L$ of Sec.~\ref{sec:weak_regime}.

Statements about the full problem use the layered structure of Sec.~\ref{sec:boundary_layers}, which we have established only formally (Remark~\ref{rem:optimal_trajectories}), so we state it as hypothesis (H): there is $t_\varepsilon = O(\ln(1/\varepsilon))$ such that on $[t_\varepsilon, T - t_\varepsilon]$ the extremal lies within $O(\varepsilon^{3/2})$ in $\bar c$ and $O(\varepsilon^{(3+\omega)/2})$ in $\bar p$ of a trajectory on $\mathcal M_\varepsilon$ whose amplitudes are $(\bar A, \bar Q) + O(\varepsilon^{1/2})$. On $[0,t_\varepsilon]\cup[T-t_\varepsilon,T]$, $\bar c = O(\varepsilon^{1/2})$ and $\bar p = O(\varepsilon^{(1+\omega)/2})$. Since the layers \eqref{eq:layer_solutions} decay at rate $\Lambda$, a choice $t_\varepsilon = K\Lambda^{-1}\ln(1/\varepsilon)$ with $K$ large is consistent with (H).

\begin{theorem}[Second variation near the organizing center]
\label{thm:second_order}
Let the hypotheses of Theorem~\ref{thm:reduction} hold with $\omega \in \{0,1,2\}$, and let $(\bar c, \bar p, \bar u)$ be an extremal whose leading-order form is given by that theorem. For $\omega>0$, assume also (H), and write $\tilde u^\alpha(t) = U^\alpha(\tau)$ with $\tau = \varepsilon^{\omega/2}t$. Every $\tilde u \in L^2(0,T)$ has this form, with $U \in L^2(0,T_\tau)$. Write $U^\alpha = e^\alpha v + U_\perp^\alpha$, where $e^\alpha = r^i{}_\alpha w_c^i/\sqrt{\beta_c}$ and $e^\alpha U_\perp^\alpha = 0$. The constants implied by $O(\cdot)$ depend on the extremal and on $T_\tau$, but not on $U$.

(i) In the strong regime, $\omega=0$, the second variation is, for all $\tilde u \in L^2(0,T)$,
\begin{equation}
    \delta^2 J = \mathcal{Q}_0[\tilde u] + O(\varepsilon^{1/2})\|\tilde u\|^2_{L^2},
    \qquad
    \mathcal{Q}_0[\tilde u] = \chi\|\tilde c(T)\|^2 + \int_0^T\big[\zeta\|\tilde c\|^2 + \|\tilde u\|^2\big]dt ,
    \label{eq:2var_strong}
\end{equation}
with $\tilde c$ solving \eqref{eq:linearized_state} with coefficients at the organizing center. Since $\mathcal{Q}_0[\tilde u] \ge \|\tilde u\|^2_{L^2}$, $\delta^2 J$ is coercive for sufficiently small $\varepsilon$.

(ii) In the intermediate regime, $\omega=1$, the second variation is
\begin{equation}
    \delta^2 J = \varepsilon^{-1/2}\Big\{\mathcal{Q}_1[v] + \|U_\perp\|^2_{L^2} + O(\varepsilon^{1/4})\|U\|^2_{L^2}\Big\},
    \label{eq:2var_int}
\end{equation}
with $\Theta = 0$ and $H_{AA} = \tilde\zeta$. Since $\mathcal{Q}_1[v] \ge \|v\|^2_{L^2}$, the leading-order form is coercive.

(iii) In the weak regime, $\omega=2$, the second variation is
\begin{equation}
    \delta^2 J = \varepsilon^{-1}\Big\{\mathcal{Q}_2[v] + \|U_\perp\|^2_{L^2} + O\big(\varepsilon^{1/2}\ln(1/\varepsilon)\big)\|U\|^2_{L^2}\Big\},
    \label{eq:2var_weak_full}
\end{equation}
with $\Theta = f'(\bar A)$ and $H_{AA} = \tilde\zeta + 6\nu\bar A\bar Q$. Let $y$ solve the Jacobi equation
\begin{equation}
\begin{gathered}
    y'' = W(\tau)\,y, \qquad
    W = \beta_c\tilde\zeta + f'(\bar A)^2 + f(\bar A)\,f''(\bar A),\\
    y(T_\tau) = 1,\qquad y'(T_\tau) = f'(\bar A(T_\tau)) - \beta_c\tilde\chi .
\end{gathered}
    \label{eq:jacobi}
\end{equation}
Then (a) if $y \neq 0$ on $[0, T_\tau]$, $\mathcal{Q}_2[v] \ge \gamma_0\,\|v\|^2_{L^2}$ for some $\gamma_0>0$, and (b) if $y(\tau_c) = 0$ for some $\tau_c \in (0, T_\tau)$, $\mathcal{Q}_2[v] < 0$ for some $v$.

In all cases the unfolding parameters enter only through the extremal and, for $\omega=2$, through $f'(\bar A)$ and $f(\bar A)$ in $W$.
\end{theorem}

\begin{proof}
Consider first $\omega > 0$, and write $\tilde c^i = \tilde a_c v_c^i + \sum_s \tilde a_s v_s^i$. The stable amplitudes have two forcings: the control, and the critical amplitude through the quadratic term evaluated along $\bar c = \varepsilon^{1/2}\bar A v_c + O(\varepsilon)$. Split $\tilde a_s = \tilde a_s^U + \tilde a_s^\xi$ accordingly.

The control-driven part,
\begin{equation*}
\tilde a_s^U(t) = \int_0^t e^{\lambda_s(t-t')}\,\overline{w_s}^{\,i} r^i{}_\alpha\tilde u^\alpha(t')\,dt',
\end{equation*}
is not approximated. Since $\mathrm{Re}\,\lambda_s \le -\Lambda$, Young's inequality and the Cauchy--Schwarz inequality give, uniformly in $\varepsilon$,
\begin{equation}
  \|\tilde a_s^U\|_{L^2(0,T)} \le C\|\tilde u\|_{L^2(0,T)}, \qquad
  |\tilde a_s^U(T)| \le C\|\tilde u\|_{L^2(0,T)}, \qquad
  \|\tilde u\|_{L^2(0,T)} = \varepsilon^{-\omega/4}\|U\|_{L^2(0,T_\tau)} .
  \label{eq:aU_bounds}
\end{equation}

The critical amplitude varies on the slow scale. Write $\hat a_c = \varepsilon^{\omega/2}\tilde a_c$. Then $\hat a_c(0) = 0$, and by \eqref{eq:critical_amplitude} below and Gronwall's inequality on the fixed interval $[0,T_\tau]$,
\begin{equation*}
    \|\hat a_c\|_{L^\infty(0,T_\tau)} + \|\hat a_c'\|_{L^1(0,T_\tau)} \le C\|U\|_{L^2(0,T_\tau)} .
\end{equation*}
The part of the stable amplitudes that $\hat a_c$ drives is therefore slaved,
\begin{equation}
    \tilde a_s^\xi = -\frac{\varepsilon^{1/2}}{\lambda_s}\,\bar A\,\eta_{sc}\,\tilde a_c + \text{h.o.t.},
    \label{eq:stable_slaving}
\end{equation}
in the following integrated sense. With $h = \varepsilon^{\omega/2}$ and $F = \bar A\hat a_c$, the leading stable filter in slow time is $K_hF(\tau) = h^{-1}\int_0^\tau e^{\lambda_s(\tau-s)/h}F(s)\,ds$, and $\tilde a_s^\xi = \varepsilon^{1/2}\eta_{sc}\,h^{-1}K_hF$ at leading order. Since $F(0) = 0$, integration by parts and Young's inequality give
\begin{equation*}
    K_hF + \frac{F}{\lambda_s} = \frac{1}{\lambda_s}\int_0^\tau e^{\lambda_s(\tau-s)/h}F'(s)\,ds,
    \qquad
    \Big\|K_hF + \frac{F}{\lambda_s}\Big\|_{L^1(0,T_\tau)} \le C h\,\|F'\|_{L^1(0,T_\tau)} \le C\varepsilon^{\omega/2}\|U\|_{L^2(0,T_\tau)},
\end{equation*}
where the last step uses the bound on $\hat a_c$ above. The slaving error is thus bounded without differentiating $U$, so neither part requires $U$ to vary slowly, and there is no initial layer. The slaved part is the tangent to the center-manifold graph and has size $\varepsilon^{(1-\omega)/2}\hat a_c$; perturbed coefficients and the endpoint intervals are handled using (H), as below.

Next, project \eqref{eq:linearized_state} onto $w_c^i$ and substitute \eqref{eq:stable_slaving}. The $O(\varepsilon^{1/2})$ part of the coefficient of $\tilde a_c$ vanishes by Assumption~\ref{assump:cusp}, and $h_1 = O(\varepsilon^{3/2})$ contributes only at $O(\varepsilon^{3/2})$. At $O(\varepsilon)$, the direct coefficient $w_c^m r^m{}_j(\bar c,\bar u,h)\,v_c^j$ contributes $\mu\tilde h_2 + \big(\tfrac12 w_c^m r^m{}_{jkl}v_c^jv_c^kv_c^l - \tfrac12\sum_s\eta_{cs}\eta_{sc}/\lambda_s\big)\bar A^2$, and the slaved part \eqref{eq:stable_slaving} contributes $-\sum_s\eta_{cs}\eta_{sc}\bar A^2/\lambda_s$. Their sum is $\varepsilon f'(\bar A)$, with $f'(\bar A) = \mu\tilde h_2 + 3\nu\bar A^2$ by \eqref{eq:nu}. With $\hat a_c = \varepsilon^{\omega/2}\tilde a_c$ and $w_c^i r^i{}_\alpha U^\alpha = \sqrt{\beta_c}\,v$,
\begin{equation}
    \hat a_c' = \varepsilon^{(2-\omega)/2} f'(\bar A)\,\hat a_c + \sqrt{\beta_c}\,v + O(\varepsilon^{1/2}), \qquad \hat a_c(0) = 0 .
    \label{eq:critical_amplitude}
\end{equation}
At leading order this is \eqref{eq:xi_lin} with $\xi = \hat a_c$. The remainder in \eqref{eq:critical_amplitude} has two sources: the forcing $\varepsilon^{1/2}\bar A\,\eta_{cs}\tilde a_s^U$, which by \eqref{eq:aU_bounds} is $O(\varepsilon^{1/2})\|U\|_{L^2}$ in $L^2(0,T_\tau)$, and the deviations of the coefficients from their outer values that (H) allows. In the layers these deviations are $O(\varepsilon^{1/2})$ over a time $t_\varepsilon$. Hence $\|\hat a_c - \xi\|_{L^\infty} \le C\varepsilon^{1/2}\ln(1/\varepsilon)\|U\|_{L^2}$.

We now substitute into \eqref{eq:second_variation}, using $dt = \varepsilon^{-\omega/2}d\tau$, $\chi = \varepsilon^{\omega/2}\tilde\chi$, $\zeta = \varepsilon^\omega\tilde\zeta$ and, on $\mathcal M_\varepsilon$, $\bar p^m = \varepsilon^{(1+\omega)/2}\big(\bar Q w_c^m + \varepsilon^{1/2}\sum_s b_s\overline{w_s}^{\,m} + \dots\big)$, where the dots denote critical and higher-order terms that do not contribute below. The control cost $\int\|\tilde u\|^2dt$, the critical tracking term $\int\zeta\tilde a_c^2\,dt$ and the terminal term $\chi\tilde a_c(T)^2$ are each $\varepsilon^{-\omega/2}$ times their slow-time counterparts. Using $\|U\|^2 = v^2 + \|U_\perp\|^2$, they give $\varepsilon^{-\omega/2}\{\tilde\chi\xi(T_\tau)^2 + \int(v^2 + \|U_\perp\|^2 + \tilde\zeta\xi^2)\}$. The leading curvature term $\bar Q\,w_c^m r^m{}_{ij}v_c^iv_c^j\,\tilde a_c^2$ vanishes by Assumption~\ref{assump:cusp}. Three curvature terms of order $\varepsilon^{(2+\omega)/2}\tilde a_c^2$ remain. Evaluating $r^m{}_{ij}$ along $\bar c$ gives the coefficient $\bar Q\bar A\,w_c^m r^m{}_{ijk}v_c^iv_c^jv_c^k$. The cross term with the slaved part \eqref{eq:stable_slaving} gives $-2\bar Q\bar A\sum_s\eta_{cs}\eta_{sc}/\lambda_s$, and the stable part of the costate correction gives $\sum_s b_s\eta_{sc}$. For $\omega = 2$, $b_s = -\eta_{cs}\bar A\bar Q/\lambda_s$ by \eqref{eq:as_weak}, and by \eqref{eq:nu} the sum is $6\nu\bar A\bar Q = f''(\bar A)\bar Q$, the $A$-curvature of the reduced Hamiltonian. For $\omega = 1$ they enter only at relative order $\varepsilon^{1/2}$ and need not be evaluated. After integration they contribute $\varepsilon^{-\omega/2}\cdot\varepsilon^{(2-\omega)/2}\int f''(\bar A)\bar Q\,\xi^2$. The remaining terms are bounded using \eqref{eq:aU_bounds} together with $\|\tilde a_c\|_{L^\infty} \le C\varepsilon^{-\omega/2}\|U\|_{L^2}$ and $\|\tilde a_c\|_{L^2(0,T)} \le C\varepsilon^{-3\omega/4}\|U\|_{L^2}$. Relative to the leading factor $\varepsilon^{-\omega/2}\|U\|^2_{L^2}$, the stable-mode tracking terms and the terminal stable term are $O(\varepsilon^{\omega/2})$, and all curvature and mixed terms containing $\tilde a_s^U$ or $\tilde u$ are $O(\varepsilon^{1/2})$. The mixed critical/stable tracking terms $2\zeta\,\tilde a_c\tilde a_s\,v_c^iv_s^i$, present when $v_c$ is not orthogonal to the stable eigenvectors, are likewise $O(\varepsilon^{1/2})$ relative to the leading factor. The terminal cross term $2\chi\,\tilde a_c(T)\,\tilde a_s^U(T)\,v_c^i v_s^i$ is $O(\varepsilon^{\omega/4})$; it vanishes when $v_c$ is orthogonal to the stable eigenvectors. The layers and transverse errors allowed by (H) contribute $O(\varepsilon^{1/2}\ln(1/\varepsilon))$, with the logarithm coming from the layer width.
Collecting the terms, for $\omega \in \{1,2\}$,
\begin{equation}
\begin{split}
    \delta^2 J = \varepsilon^{-\omega/2}\Big\{&\tilde\chi\,\xi(T_\tau)^2 + \int_0^{T_\tau}\big[v^2 + \|U_\perp\|^2 + \big(\tilde\zeta + \varepsilon^{(2-\omega)/2}f''(\bar A)\bar Q\big)\xi^2\big]d\tau \\
    &+ O\big(r_\omega(\varepsilon)\big)\|U\|^2_{L^2}\Big\},
\end{split}
    \label{eq:2var_general}
\end{equation}
with $r_1 = \varepsilon^{1/4}$ and $r_2 = \varepsilon^{1/2}\ln(1/\varepsilon)$. Since $\varepsilon^{(2-\omega)/2} = \theta_\omega + O(\varepsilon^{1/2})$ for $\omega\in\{1,2\}$, this gives \eqref{eq:2var_int} and \eqref{eq:2var_weak_full}. Every $\tilde u \in L^2(0,T)$ is of the form $U(\varepsilon^{\omega/2}t)$, so these expansions hold for all perturbations.

For $\omega=0$ there is no scale separation, and $T$ is fixed. The curvature terms carry $\bar p = O(\varepsilon^{1/2})$, and the coefficients of \eqref{eq:linearized_state} are those at the organizing center up to $O(\varepsilon^{1/2})$. Since $\|\tilde c\|_{L^\infty} \le C\|\tilde u\|_{L^2}$ on $[0,T]$ by Gronwall's inequality, these corrections are $O(\varepsilon^{1/2})\|\tilde u\|^2_{L^2}$ for every $\tilde u \in L^2(0,T)$, which gives \eqref{eq:2var_strong}. Coercivity follows because $\chi,\zeta \ge 0$. For $\omega=1$, $H_{AA} = \tilde\zeta > 0$, so $\mathcal Q_1[v]\ge\|v\|^2_{L^2}$.

For case (iii), the Jacobi system is the linearization of \eqref{eq:Aweak}--\eqref{eq:Qweak}:
\begin{equation}
    y' = \Theta y - \beta_c z, \qquad z' = -H_{AA}\,y - \Theta z, \qquad z(T_\tau) = \tilde\chi\,y(T_\tau).
    \label{eq:jacobi_system}
\end{equation}
Eliminating $z = (\Theta y - y')/\beta_c$ gives $y'' = (\Theta' + \Theta^2 + \beta_c H_{AA})\,y$. To evaluate the coefficient, use $\bar A' = f(\bar A) - \beta_c\bar Q$:
\begin{equation*}
    \Theta' = f''(\bar A)\,\bar A' = f''f - \beta_c f''\bar Q,
    \qquad
    \beta_c H_{AA} = \beta_c\tilde\zeta + \beta_c f''\bar Q .
\end{equation*}
The costate terms cancel, which leaves the coefficient $W$ of \eqref{eq:jacobi}. The terminal data follow from $y'(T_\tau) = \Theta(T_\tau) - \beta_c z(T_\tau)$.

For (a), $\Pi = z/y$ is $C^1$ on $[0,T_\tau]$ and satisfies the Riccati equation $-\Pi' = 2\Theta\Pi - \beta_c\Pi^2 + H_{AA}$ with $\Pi(T_\tau) = \tilde\chi$. Differentiating $\Pi\xi^2$ along \eqref{eq:xi_lin}, integrating over $[0,T_\tau]$, and adding the result to \eqref{eq:2var_weak} completes the square:
\begin{equation*}
    \mathcal{Q}_2[v] = \int_0^{T_\tau}\big(v + \sqrt{\beta_c}\,\Pi\,\xi\big)^2 d\tau .
\end{equation*}
Coercivity then follows from the Riccati sufficient condition \cite[Thm.~5.2]{maurer1981first}, of which this is the scalar case.

For (b), let $\hat v = 0$ on $[0,\tau_c)$ and $\hat v = -\sqrt{\beta_c}\,z$ on $[\tau_c, T_\tau]$. Its response is $\xi = 0$ on $[0,\tau_c]$ and $\xi = y$ on $[\tau_c,T_\tau]$. Using $(yz)' = -\beta_c z^2 - H_{AA}y^2$ and $y(\tau_c)=0$,
\begin{equation*}
    \mathcal{Q}_2[\hat v] = \tilde\chi - \big[yz\big]_{\tau_c}^{T_\tau} = 0 .
\end{equation*}
If $\mathcal{Q}_2 \ge 0$, then $\hat v$ minimizes $\mathcal{Q}_2$, and its first variation vanishes. This forces $\hat v = -\sqrt{\beta_c}\,\hat z$ for a continuous adjoint $\hat z$. Since $\hat v = 0$ on $[0,\tau_c)$, $\hat z$ vanishes there; on $[\tau_c,T_\tau]$ it solves the same equation as $z$ with the same terminal value, so $\hat z = z$. Continuity at $\tau_c$ gives $z(\tau_c)=0$, which is impossible because $y$ and $z$ cannot vanish simultaneously. Because the right endpoint is free and carries a terminal penalty, $\tau_c$ is a focal point rather than a conjugate point, and this is the corresponding Jacobi necessary condition \cite[Ch.~5]{gelfand1963calculus}, \cite{hestenes1966calculus}.
\end{proof}

\begin{theorem}[Local optimality]
\label{thm:local_optimality}
Let the hypotheses of Theorem~\ref{thm:second_order} hold, with $y$ the Jacobi field \eqref{eq:jacobi}.

(i) If $\omega = 2$ and $y(\tau_c) = 0$ for some $\tau_c \in (0,T_\tau)$, then for all sufficiently small $\varepsilon$ the extremal is not a weak local minimizer.

(ii) If $\omega \in \{0,1\}$, or $\omega = 2$ and $y \ne 0$ on $[0,T_\tau]$, there is $\gamma > 0$, independent of $\varepsilon$, such that $\delta^2 J[\tilde u] \ge \gamma\,\|\tilde u\|^2_{L^2(0,T)}$ for all $\tilde u \in L^2(0,T)$ and all sufficiently small $\varepsilon$. For each such $\varepsilon$ the extremal is a strict weak local minimizer: there are $\alpha, \rho > 0$ such that $J[u] \ge J[\bar u] + \alpha\|u - \bar u\|^2_{L^2}$ whenever $\|u - \bar u\|_{L^\infty} \le \rho$.
\end{theorem}

\begin{proof}
For (i), $\mathcal Q_2$ is continuous on $L^2$, so by Theorem~\ref{thm:second_order}(iii)(b) there is a smooth $v$ with $\mathcal Q_2[v] < 0$, and \eqref{eq:2var_weak_full} with $U = ev$ gives $\delta^2 J[\tilde u] < 0$ for $\tilde u(t) = U(\varepsilon t)$ for small $\varepsilon$. Since $J$ is $C^2$ on $L^\infty$, the function $s \mapsto J[\bar u + s\tilde u]$ then has a strict local maximum at $s = 0$, so $\bar u$ is not a weak local minimizer.

For (ii), let $\gamma_0$ be the coercivity constant of $\mathcal Q_\omega$: $\gamma_0 = 1$ for $\omega = 1$, and $\gamma_0$ is given by Theorem~\ref{thm:second_order}(iii)(a) for $\omega = 2$. Then $\mathcal Q_\omega[v] + \|U_\perp\|^2_{L^2} \ge \min(\gamma_0,1)\|U\|^2_{L^2}$. Since $\varepsilon^{-\omega/2}\|U\|^2_{L^2(0,T_\tau)} = \|\tilde u\|^2_{L^2(0,T)}$, the expansions \eqref{eq:2var_int} and \eqref{eq:2var_weak_full} give $\delta^2 J[\tilde u] \ge \tfrac12\min(\gamma_0,1)\|\tilde u\|^2_{L^2(0,T)}$ for all $\tilde u$ once $\varepsilon$ is small. For $\omega = 0$ the same bound follows from Theorem~\ref{thm:second_order}(i). Local optimality then follows from the two-norm second-order sufficient condition \cite[Thm.~3.5]{maurer1981first}.
\end{proof}

Statements (iii)(a)--(b) of Theorem~\ref{thm:second_order} are exact for the reduced problem. By \cite[Thm.~5.2]{maurer1981first} and the argument used for Theorem~\ref{thm:local_optimality}(i), the reduced extremal $\bar v = -\sqrt{\beta_c}\bar Q$ is a strict weak local minimizer if $\mathcal Q_2$ is coercive, and is not a weak local minimizer if $\mathcal Q_2[v] < 0$ for some $v$. If $y(0) = 0$ and $y \ne 0$ on $(0,T_\tau]$, then $\mathcal Q_2 \ge 0$ but $\mathcal Q_2[-\sqrt{\beta_c}z] = 0$, and the test is inconclusive.

In the weak regime, applying Theorem~\ref{thm:second_order}(iii) requires the zeros of the Jacobi field $y$, which in general must be computed. The extremal itself, however, already constrains where they can lie. Because $W$ depends only on $\bar A$, the time-translation mode of the extremal supplies a second solution of \eqref{eq:jacobi}, which can be compared with $y$.

\begin{lemma}[Time-translation Jacobi field]
\label{lem:translation}
Along any extremal of \eqref{eq:Aweak}--\eqref{eq:Qweak} (with $A^*$ constant),
\begin{equation}
    \bar A'' = \beta_c\tilde\zeta\,(\bar A - A^*) + f(\bar A)\,f'(\bar A),
    \label{eq:newton}
\end{equation}
and $\bar A'$ solves $y'' = W y$. If $\bar A$ is not constant, the zeros of $\bar A'$ are simple and isolated.
\end{lemma}

\begin{proof}
Differentiate $\bar A' = f - \beta_c\bar Q$ and substitute \eqref{eq:Qweak}:
\begin{equation*}
    \bar A'' = f'\bar A' + \beta_c\tilde\zeta(\bar A - A^*) + \beta_c f'\bar Q = f'(\bar A' + \beta_c\bar Q) + \beta_c\tilde\zeta(\bar A - A^*),
\end{equation*}
and $\bar A' + \beta_c\bar Q = f$. Differentiating \eqref{eq:newton} once more gives $(\bar A')'' = W\bar A'$, since $W$ is the $A$-derivative of the right side of \eqref{eq:newton}. A double zero of $\bar A'$ forces $\bar A'\equiv0$ by uniqueness for linear ODEs.
\end{proof}

Equation \eqref{eq:newton} is Newton's equation for motion in the potential $\tfrac{\beta_c\tilde\zeta}{2}(A - A^*)^2 + \tfrac12 f(A)^2$, and $W$ is the curvature of this potential along the extremal. This motivates the canonical change of variables of Sec.~\ref{sec:canonical}, which brings the reduced Hamiltonian to natural form with this potential.

\begin{corollary}[Turning points]
\label{cor:turning}
Let $\omega=2$ and let $\bar A$ be nonconstant. If $\bar A'(\tau_1) = \bar A'(\tau_2) = 0$ for some $0 < \tau_1 < \tau_2 \le T_\tau$, then $\mathcal{Q}_2$ takes negative values.
\end{corollary}

\begin{proof}
By Lemma~\ref{lem:translation}, $\bar A'$ is a nontrivial solution of $y''=Wy$ that vanishes at $\tau_1$ and $\tau_2$. By the Sturm separation theorem \cite[Ch.~XI]{hartman2002ode}, the solution $y$ of \eqref{eq:jacobi} vanishes somewhere in $[\tau_1,\tau_2]$. Since $y(T_\tau)=1$, this zero lies in $[\tau_1,T_\tau)\subset(0,T_\tau)$, and Theorem~\ref{thm:second_order}(iii)(b) applies.
\end{proof}

\section{Asymptotic structure of the weak-control regime}
\label{sec:nonlinear_analysis}

In this section we analyze solutions of the first-order optimality conditions \eqref{eq:Aweak}--\eqref{eq:Qweak} in the weak-control regime. As a first step, we bring the reduced Hamiltonian to natural form.

\subsection{Canonical variables}
\label{sec:canonical}

Write $\hat\mu \equiv \mu\tilde h_2$ for the linear coefficient of the reduced drift $f(A) = \sigma\tilde h_1 + \hat\mu A + \nu A^3$, and assume $\hat\mu \neq 0$. Let $L \equiv \sqrt{|\hat\mu/\nu|}$ and define
\begin{equation}
    P \equiv \frac{\beta_c Q - f(A)}{|\nu| L^3}, \qquad
    X \equiv \frac{A}{L}, \qquad
    \hat\tau \equiv \tau/T_{\tau}.
\end{equation}
Since $A' = f(A) - \beta_c Q$, the momentum $P = -A'/(|\nu| L^3)$ is proportional to the velocity. Two identities are needed. First, the reduced Hamiltonian satisfies $H_{\text{reduced}} = -(\beta_c Q - f)^2/(2\beta_c) + \tfrac{\tilde\zeta}{2}(A-A^*)^2 + f^2/(2\beta_c)$, as expanding the square shows. Second, substituting $A = LX$ gives $f(LX) = \nu L^3\,g(X)$, with $g(X) \equiv X^3 + aX + b$ and $a$, $b$ as below.
In these variables,
\begin{equation}
    \mathcal{H} \equiv \frac{H_{\text{reduced}}}{J} = -\frac{1}{2}P^2 + V(X),
    \qquad
    V(X) = \frac{\kappa}{2}(X-X^*)^2 + \frac{1}{2}g(X)^2 ,
    \label{eq:canonical_potential}
\end{equation}
with the rescaled parameters
\begin{equation}
\begin{gathered}
    J = \frac{|\hat\mu|^3}{\beta_c |\nu|}, \qquad
    S = T_\tau |\hat\mu|,
    \qquad
    \kappa = \frac{\tilde{\zeta}\,\beta_c}{\hat\mu^2},
    \qquad
    X^* = \frac{A^*}{L},
    \\
    a = \operatorname{sgn}(\hat\mu\nu),
    \qquad
    b = \frac{\sigma\tilde h_1}{\nu L^3} = \operatorname{sgn}(\nu)\,\frac{\sigma\tilde h_1 |\nu|^{1/2}}{|\hat\mu|^{3/2}} .
\end{gathered}
    \label{eq:canonical_params}
\end{equation}
The change of variables multiplies the symplectic form by the constant $\beta_c/(|\nu|L^4)$ and the Hamiltonian by $1/J$. Since $J\beta_c/(|\nu|L^4) = |\hat\mu|$, the pair $(X,P)$ is canonical, and evolution in $\hat\tau$ is generated by $S\mathcal H$. With overdots denoting $d/d\hat\tau$,
\begin{align}
    \dot{X} &= S\frac{\partial \mathcal{H}}{\partial P} = -S P, \\
    \dot{P} &= -S\frac{\partial \mathcal{H}}{\partial X} = -\kappa S(X-X^*) - S\,g(X)\,g'(X).
\end{align}
The terminal condition $Q(T_\tau) = \tilde\chi(A(T_\tau) - A^*)$ becomes, at $\hat\tau=1$,
\begin{equation}
    P = \mu_b (X-X^*) - \operatorname{sgn}(\nu)\,g(X),
    \qquad \mu_b = \frac{\beta_c\tilde{\chi}}{|\hat\mu|} .
    \label{eq:canonical_terminal}
\end{equation}

The boundary-value problem is therefore parameterized by $X_0$, $\kappa$, $\mu_b$, $S$, $X^*$ and $b$, together with the binary choices $a = \pm 1$ and $\operatorname{sgn}(\nu)$. Only $(\kappa, X^*, a, b)$ affect the phase portrait. The terminal condition is set by $\mu_b$ and $\operatorname{sgn}(\nu)$, while $S$ and $X_0$ select the solution of the boundary-value problem. The flow, potential and terminal condition are invariant under $(X, P, X_0, X^*, b) \mapsto (-X, -P, -X_0, -X^*, -b)$, so $b \ge 0$ may be assumed without loss of generality.

The second-order quantities of Sec.~\ref{sec:second_order} take a simple form in these variables. The Jacobi coefficient of \eqref{eq:jacobi} is $W = \hat\mu^2\,V_{XX}(\bar X)$, so $W>0$ near saddles of the canonical flow and $W<0$ near centers. Since $P = -\bar A'/(|\nu|L^3)$, the turning points of Corollary~\ref{cor:turning} are crossings of the axis $P = 0$, and the corollary states that an extremal which is locally optimal at leading order crosses $P = 0$ at most once in $0 < \hat\tau \le 1$.

The free part of the potential, $\tfrac12 g^2$, vanishes exactly at the equilibria of the uncontrolled amplitude equation \eqref{eq:cusp_amplitude}. For $a = -1$ there are three such zeros when $|b| < b_c \equiv 2/(3\sqrt3)$, and one otherwise. The value $|b| = b_c$ is the fold curve $\Delta = 0$ of Eq.~\eqref{eq:cusp_discriminant}. For $a = +1$ there is always exactly one zero. On the ray $\tilde h_2 = 0$, where $\hat\mu = 0$, and provided $\sigma\tilde h_1 \neq 0$, the same form holds with $a = 0$, $b = \pm 1$ and $L = |\sigma\tilde h_1/\nu|^{1/3}$, once $|\hat\mu|$ is replaced by $|\nu|L^2$ in \eqref{eq:canonical_params} and in $\mu_b$ (the two agree when $\hat\mu \neq 0$). The organizing center $\tilde h_1 = \tilde h_2 = 0$ itself is excluded from this normalization.

\subsection{Bifurcations}
\label{sec:bifurcation}

We first restrict to the pitchfork ray $b = 0$, $a = -1$, i.e.\ $\tilde h_1 = 0$ with $\hat\mu\nu < 0$, where
\begin{equation}
    V(X;\kappa,X^*)=\frac{\kappa}{2}(X-X^*)^2+\frac{1}{2}(X-X^3)^2,
\end{equation}
and then follow the codimension-two points of this diagram off the ray.

Equilibria are the points $(X_e, 0)$ with $V_X(X_e) = 0$. Since the flow is Hamiltonian with one degree of freedom, a nondegenerate equilibrium is of one of two types: a saddle if $V_{XX}(X_e) > 0$ (a local minimum of $V$), or a center if $V_{XX}(X_e) < 0$ (a local maximum of $V$). Two codimension-one bifurcations organize the $(\kappa, X^*)$ plane (Fig.~\ref{fig:local_bifs}). At a saddle-node (green curves), $V_X = V_{XX} = 0$ and a saddle and a center coalesce. At a heteroclinic connection (blue curves), two neighboring saddles separated by a single center lie at the same level of $V$. Crossing this curve leaves the equilibria unchanged but switches the saddle to which the homoclinic loop around the center is attached. The curves meet at three types of codimension-two point (markers): cusps, where $V_X = V_{XX} = V_{XXX} = 0$; double saddle-nodes, where $V_X = V_{XX} = 0$ at two distinct points; and intersections of a heteroclinic curve with a saddle-node curve.

\begin{figure}[!t]
    \centering
    \includegraphics{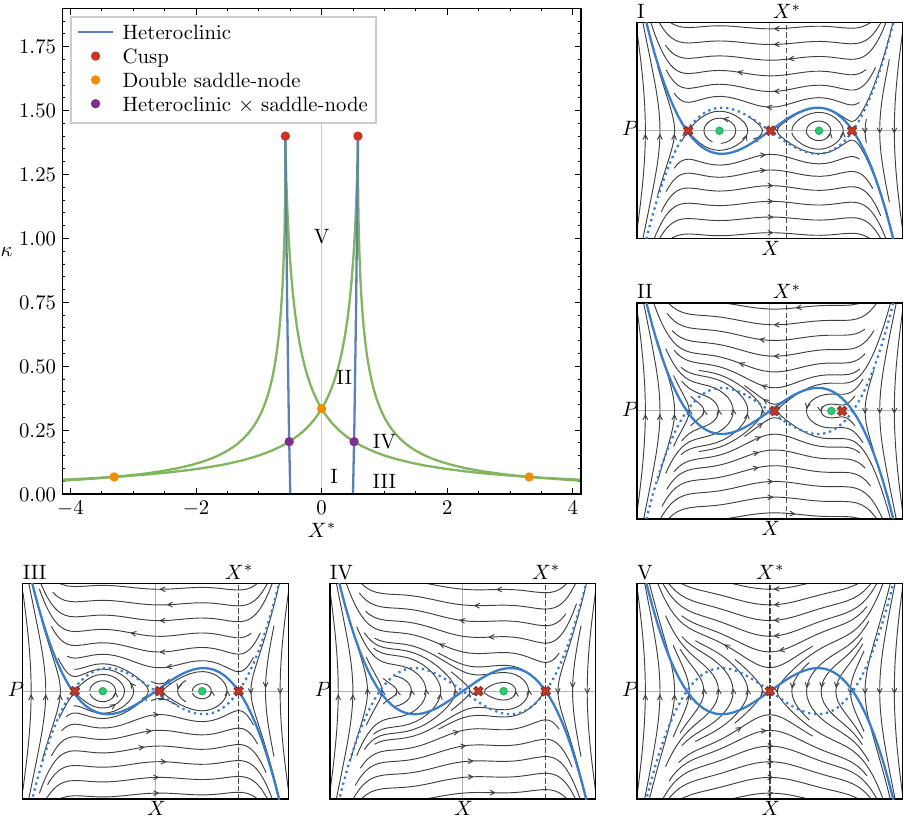}
    \caption{Bifurcation diagram of the canonical Hamiltonian flow on the pitchfork ray $b = 0$, $a = -1$, in the $(X^*, \kappa)$ parameter plane (upper left), with representative phase portraits in each labeled region (I--V). Green curves mark saddle-node bifurcations and blue curves heteroclinic connections; markers are the codimension-two points. In each phase portrait, green dots are centers and red crosses are saddles. Blue curves are separatrices, and the dashed vertical line marks the target $X^*$.}
    \label{fig:local_bifs}
\end{figure}

The number and type of equilibria are constant in each region I--V of Fig.~\ref{fig:local_bifs}. For large $\kappa$ the potential is convex and the flow has a single saddle (region V). The saddle-node curve bounding region V can be interpreted as the boundary between the weak-control and intermediate-control regimes: above it, the flow is topologically equivalent to the linear saddle of the intermediate regime. For the control problem, saddles are the candidate turnpikes of Sec.~\ref{sec:boundary_layer}, while orbits near a center oscillate in $X$; by Corollary~\ref{cor:turning}, such orbits are not locally optimal once they have two turning points in $(0,T_\tau]$.

For $b \neq 0$ the reflection symmetry is broken, and each codimension-two point traces a curve in $(X^*, b, \kappa)$ (Fig.~\ref{fig:codim2_tracks}). By the symmetry $(X^*, b) \mapsto (-X^*, -b)$ it suffices to consider $b > 0$. The configuration of Fig.~\ref{fig:local_bifs} persists for small $b$, and the codimension-two points disappear in turn as $b$ increases. At $b = b_\infty$, an outer double saddle-node reaches $\kappa = 0$ and escapes to $|X^*| \to \infty$. At $b = b_\star$, a cusp crosses the saddle-node curve of a second degenerate point; at this codimension-three point the remaining double saddle-nodes annihilate, and one heteroclinic--saddle-node intersection ends in the cusp. At $b = b_c$, a cusp and the remaining intersection reach $\kappa = 0$. For $b > b_c$, where the uncontrolled amplitude equation has a single equilibrium, only one cusp remains.

\begin{figure}[!t]
    \centering
    \includegraphics{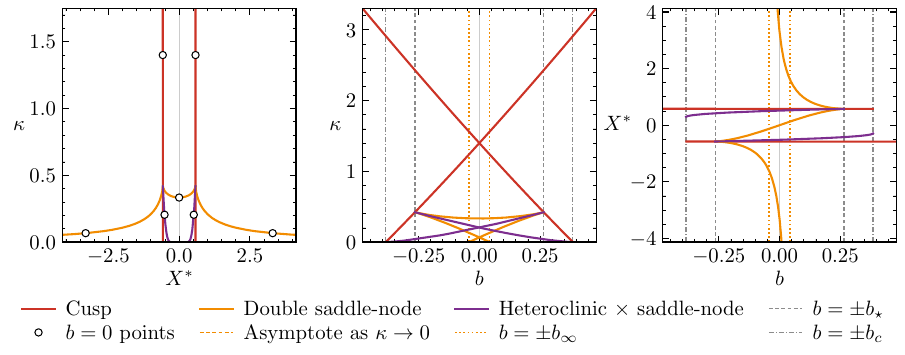}
    \caption{Continuation in $b$ of the codimension-two points of Fig.~\ref{fig:local_bifs}, projected onto the $(X^*,\kappa)$ (left), $(b,\kappa)$ (center) and $(b,X^*)$ (right) planes. Open circles mark the points at $b = 0$, and dashed orange curves are the leading-order asymptotes of the double saddle-node curve as $\kappa \to 0$. The computation is outlined in Appendix~\ref{app:bifurcation_numerics}.}
    \label{fig:codim2_tracks}
\end{figure}

\subsection{Boundary layer analysis}
\label{sec:boundary_layer}

To our knowledge, no general closed-form solution is available for the Hamiltonian boundary-value problem of the weak regime. In the small-$\kappa$, long-horizon limit, however, matched asymptotic expansions yield accurate approximations. This is the regime of turnpike theory: for long horizons, optimal trajectories consist of a transient initial arc, a long arc exponentially close to the steady state of an associated static problem, and a transient terminal arc \cite{porretta2013long,trelat2015turnpike}. Geometrically, the turnpike is a hyperbolic saddle of the optimality Hamiltonian, and the transient arcs lie on its stable and unstable manifolds \cite{sakamoto2021geometric}. Here the local minima of $V$, which are saddles of the canonical flow (Sec.~\ref{sec:bifurcation}), serve as candidate turnpikes, and an extremal may visit several of them in sequence \cite{porretta2016remarks}. In the canonical phase plane such an extremal shadows a chain of heteroclinic connections between these saddles, and the positions of the transitions are fixed by exponentially weak interactions between neighboring layers, as in Lin's method \cite{lin1990using,homburg2010homoclinic} and for metastable patterns in scalar reaction--diffusion equations \cite{carr1989metastable,fusco1989slow}. We follow Rubinstein, Sternberg and Keller \cite{rubinstein1993front}, who treat the analogous spatial problem for a tristable potential with weakly perturbed wells. Their matched-asymptotic construction requires only the decay rates and tail amplitudes of each layer, not the layer profile itself. The same structure arises in the large-deviation theory of noisy multistable systems \cite{freidlin1970small, grafke2015instanton}.

Substituting $\ddot X = -S \dot P$ yields
\begin{equation}
    S^{-2} \ddot X(\hat\tau) = V_X(X(\hat\tau)),
    \label{eqn:schrodinger}
\end{equation}
subject to $X(0) = X_0$ and, from \eqref{eq:canonical_terminal},
\begin{equation}
    -S^{-1} \dot X(1) = \hat{\chi} (X(1) - X^*) - \operatorname{sgn}(\nu)\, g\bigl(X(1)\bigr),
\end{equation}
with $\hat\chi \equiv \mu_b$. We take $\hat\chi$ large, so that, provided the terminal momentum $P$ stays bounded, the terminal state is $X_f = X^*$ at leading order.

We split the potential as $V(X)=V_0(X)+\kappa V_1(X)$, with $V_0 = \tfrac{1}{2}g^2$ and $V_1 = \tfrac{1}{2}(X - X^*)^2$, and treat $\kappa$ as small. The minima of $V_0$ are the real zeros of $g$, and all of them lie at the level $V_0 = 0$ for every $b$. We restrict to $a = -1$ and $|b| < b_c$, where $g$ has three simple zeros
\begin{equation}
    r_1 < r_2 < r_3, \qquad
    r_k = \frac{2}{\sqrt3}\cos\!\Big(\tfrac13\arccos\!\big(-\tfrac{3\sqrt3}{2}\,b\big) - \tfrac{2\pi (3-k)}{3}\Big),
    \label{eq:roots}
\end{equation}
which reduce to $\{-1, 0, 1\}$ at $b = 0$. For large $S$, the leading-order outer solution is piecewise constant, with values in $\{r_1, r_2, r_3\}$ between internal layers.

Near a transition centered at $\hat\tau_i$, with the stretched variable $\varsigma = S(\hat\tau-\hat\tau_i)$, the layer solves the heteroclinic problem
\begin{equation}
    X''(\varsigma) = V_0'(X), \qquad X(-\infty) = X_i^-, \qquad X(+\infty) = X_i^+,
    \label{eq:inner_leading}
\end{equation}
where $X_i^\pm$ are adjacent zeros of $g$. On the zero level of $V_0$, energy conservation gives $X' = \pm g(X)$, with the sign chosen so that the layer runs from $X_i^-$ to $X_i^+$. Since the zeros are simple, partial fractions give the profile implicitly,
\begin{equation}
    \Phi(X(\varsigma)) = \Phi(X_m) \pm \varsigma, \qquad \Phi(X) \equiv \sum_{k=1}^{3} \frac{\ln|X - r_k|}{g'(r_k)}, \qquad g'(r_k) = \prod_{l \neq k}(r_k - r_l),
    \label{eq:inner_implicit}
\end{equation}
where $X_m = X(0)$ fixes the translation. Following \cite{rubinstein1993front}, internal layers are normalized at the midpoint $X_m = (X_i^- + X_i^+)/2$. On the pitchfork ray, \eqref{eq:inner_implicit} can be inverted:
\begin{equation}
    X(\varsigma) = \frac{\pm 1}{\sqrt{1 + Ce^{-2\varsigma}}}, \qquad
    P(\varsigma) = \frac{\mp Ce^{-2\varsigma}}{(1 + Ce^{-2\varsigma})^{3/2}},
\end{equation}
with $C = 3$ for the midpoint normalization. These are the layers leaving $X = 0$; layers entering $X = 0$ from $\pm 1$ follow by $\varsigma \mapsto -\varsigma$ (with the sign of $P$ reversed). For $b \neq 0$, \eqref{eq:inner_implicit} cannot in general be inverted in closed form.

Only the exponential tails of each layer enter the interaction law. As $X \to r_j$, Eq.~\eqref{eq:inner_implicit} gives
\begin{equation}
\begin{gathered}
    |X(\varsigma) - r_j| \sim \Gamma_j(X_m)\, e^{-m_j|\varsigma|},
    \qquad m_j \equiv |g'(r_j)|,
    \\
    \Gamma_j(X_m) = |X_m - r_j| \prod_{k \neq j}\left(\frac{|X_m - r_k|}{|r_j - r_k|}\right)^{g'(r_j)/g'(r_k)} .
\end{gathered}
    \label{eq:tail_amplitude}
\end{equation}
We write the tails of layer $i$ as $X - X_i^- \sim \gamma_i^L e^{m_{L,i}\varsigma}$ as $\varsigma\to-\infty$ and $X - X_i^+ \sim -\gamma_i^R e^{-m_{R,i}\varsigma}$ as $\varsigma\to+\infty$, each with the decay rate of the corresponding saddle. The signed amplitudes are
\begin{equation}
    \gamma_i^L = \operatorname{sgn}(X_m - X_i^-)\,\Gamma_{X_i^-}(X_m),
    \qquad
    \gamma_i^R = \operatorname{sgn}(X_i^+ - X_m)\,\Gamma_{X_i^+}(X_m).
    \label{eq:signed_tails}
\end{equation}
 Initial and terminal layers are given by the same formulas, with $X_m$ replaced by $X_0$ or $X_f$.

\begin{figure}
    \centering
    \includegraphics{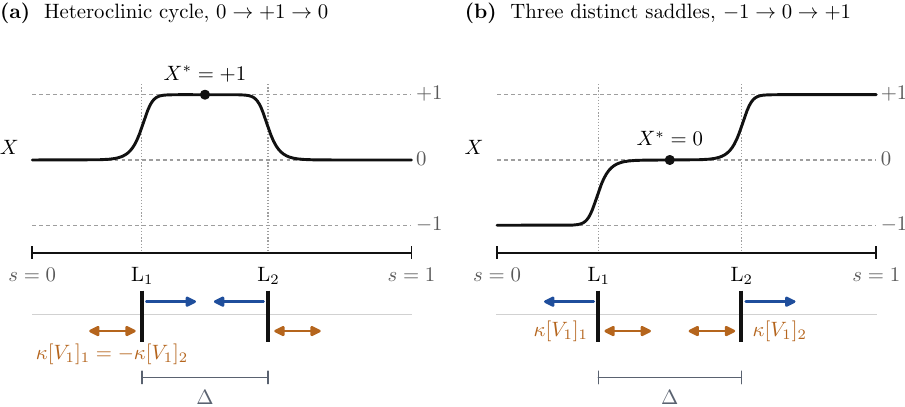}
    \caption{Schematic of interactions between neighboring internal layers. a) Fronts forming a heteroclinic cycle attract. b) Fronts joining three distinct saddles mutually repel. The drift term $\kappa [V_1]$ may move a front in either direction, depending on how close each saddle is to the target state. }
    \label{fig:front_diagram}
\end{figure}

The leading-order problem \eqref{eq:inner_leading} leaves each layer position $\hat\tau_i$ undetermined. The positions are fixed at next order by the projection argument of \cite[Sec.~2]{rubinstein1993front}: a composite approximation of neighboring layers is substituted into \eqref{eqn:schrodinger}, and the residual is projected onto the translation mode of each layer \cite{carr1989metastable,fusco1989slow}. For internal layer $i$ this gives
\begin{equation}
    \kappa [V_1]_i + 2\left(m_{L,i}^2\, \gamma_i^L \gamma_{i-1}^R\, e^{-m_{L,i} S (\hat\tau_i - \hat\tau_{i-1})} - m_{R,i}^2\, \gamma_i^R \gamma_{i+1}^L\, e^{-m_{R,i} S (\hat\tau_{i+1} - \hat\tau_i)}\right) = 0,
    \label{eq:force_balance}
\end{equation}
with decay rates $m_{L,i} = |g'(X_i^-)|$ and $m_{R,i} = |g'(X_i^+)|$, tail amplitudes from \eqref{eq:tail_amplitude}--\eqref{eq:signed_tails} with the roots \eqref{eq:roots}, and tracking jumps $[V_1]_i = V_1(X_i^+) - V_1(X_i^-)$. In \cite{rubinstein1993front} the variable analogous to $\hat\tau$ is spatial and the balance yields a velocity law. Here $\hat\tau$ is the time, and \eqref{eq:force_balance} is the solvability condition for the correction. All coefficients in \eqref{eq:force_balance} are explicit for every $|b| < b_c$. The construction requires $m S \Delta \gg 1$ for every gap $\Delta$ between layers, and it breaks down as $|b| \to b_c$, where two zeros of $g$ merge and their decay rates vanish. A schematic of the interactions is given in Fig.~\ref{fig:front_diagram}.

Eq.~\eqref{eq:force_balance} is the leading-order form of the jump conditions of Lin's method \cite{lin1990using,homburg2010homoclinic}: $\kappa[V_1]_i = \int X_i'\,\kappa V_1'(X_i)\,d\varsigma$ is the Melnikov integral of the perturbation along layer $i$, and the exponential terms are the tail interactions with its neighbors. We have not verified the hypotheses of that method with the boundary conditions at $\hat\tau = 0, 1$, and use \eqref{eq:force_balance} as a formal leading-order law.

\begin{remark}
If $X_0$ or $X_f$ is not a zero of $g$, an initial or terminal layer of the same form is included, with $X_m$ replaced by $X_0$ or $X_f$, and it interacts with the neighboring internal layer through \eqref{eq:force_balance}. If $X_0$ or $X_f$ is a zero of $g$, the interaction with the boundary is captured by a method-of-images construction, with a reflected layer across the endpoint.
\end{remark}

\begin{figure}[!t]
    \centering
    \includegraphics{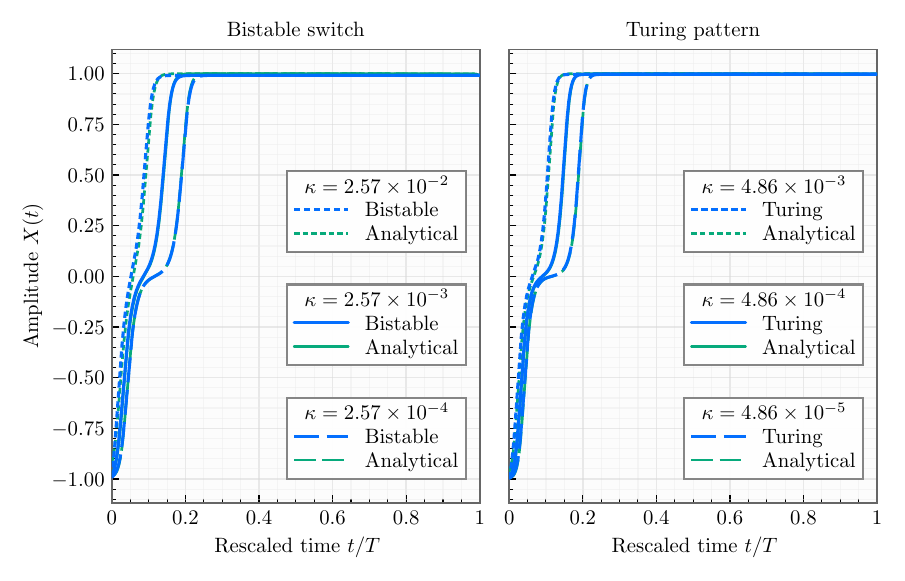}
    \caption{Boundary-layer validation on the pitchfork ray ($b = 0$): comparison between the projected amplitude trajectory from the full optimal-control solution and the matched-asymptotic profile over a range of $\kappa$, for the bistable switch (left) and the Turing pattern (right). For the bistable switch, $S \approx 59.22$ and $\hat{\chi} \approx 760$; for the Turing pattern, $S \approx 70.29$ and $\hat{\chi} \approx 45.5$. In both cases $X_0 = -X^*$, with $X^* = 0.9919$ for the bistable switch and $X^* = 0.9971$ for the Turing pattern, and $\hat{\chi}$ is large enough that $X_f = X^*$ at leading order.}
    \label{fig:boundary_layer_comparison}
\end{figure}

For $n$ internal layers, \eqref{eq:force_balance} is a system of $n$ transcendental equations for the ordered positions $\hat\tau_1 < \cdots < \hat\tau_n$, coupled only between nearest neighbors, with the endpoints $\hat\tau=0,1$ acting as neighbors of the first and last layers. When all transitions move toward the target ($[V_1]<0$), the layers cluster near $\hat\tau = 0$, and for long horizons the terminal layer contributes only an exponentially small correction to the last internal layer. The system can then be solved by a backward sweep. Define the gaps $\Delta_i \equiv \hat\tau_i - \hat\tau_{i-1}$, with $\hat\tau_0 \equiv 0$ the initial layer and $\gamma_0^R \equiv \gamma_{\mathrm{init}}^R$ its right-tail amplitude. For the last internal layer,
\begin{equation}
    \Delta_n = \frac{1}{m_{L,n} S} \ln\left(\frac{2m_{L,n}^2\, \gamma_n^L \gamma_{n-1}^R}{-\kappa [V_1]_n}\right),
    \label{eq:delta_n}
\end{equation}
and for $i = n-1, \ldots, 1$,
\begin{equation}
     \Delta_i = \frac{1}{m_{L,i} S} \ln\left(\frac{2m_{L,i}^2\, \gamma_i^L \gamma_{i-1}^R}{-\kappa [V_1]_i + 2m_{R,i}^2\, \gamma_i^R \gamma_{i+1}^L\, e^{-m_{R,i} S \Delta_{i+1}}}\right).
    \label{eq:recursion_general}
\end{equation}
Fig.~\ref{fig:boundary_layer_comparison} compares the resulting profile with the numerically computed optimal trajectories of the bistable-switch and Turing-pattern testbeds, projected onto the slow manifold and transformed to the canonical variable $X$.

The limit $S \to \infty$ also has a variational interpretation. When $v = (A' - f(A))/\sqrt{\beta_c}$ is eliminated from the reduced cost $\mathcal L$ of Sec.~\ref{sec:weak_regime}, the cross term $-\beta_c^{-1}\int f(A)A'\,d\tau$ integrates to a function of the endpoints. The reduced cost then takes the Modica--Mortola form
\begin{equation}
    \mathcal{L} = T_\tau J \int_0^1 \Big[\frac{1}{2S^2}\dot X^2 + V(X)\Big]\, d\hat\tau + \mathcal{B}\big(X(0), X(1)\big),
    \label{eq:modica_mortola}
\end{equation}
where $\mathcal{B}$ collects the endpoint terms. This is the multiwell energy of phase-transition theory, which motivates the piecewise-constant outer solution.

The results of this section give the approximate form of extremal trajectories. Sec.~\ref{sec:second_order} partly decides which of them are local minimizers, since a direction reversal of a layered extremal is a zero of $\bar A'$. By Corollary~\ref{cor:turning}, an exact reduced extremal that reverses direction two or more times is not a local minimizer of the reduced problem; for a formally assembled layer sequence, this conclusion applies to the exact extremal it approximates, when one exists. We leave extremals with at most one reversal undetermined. Identifying the global optimum among the local minimizers requires comparing their costs, which is straightforward numerically for a given sequence of fronts.

\section{Conclusion}
\label{sec:conclusion}

Near a cusp point, the leading-order extremals of a broad class of nonlinear optimal control problems take one of three forms, selected by the scaling of the tracking penalties relative to the distance from criticality (Table~\ref{scal_tab}). Under strong control, the leading-order problem is linear--quadratic tracking for the linearized system on the full state space. Under intermediate control, extremals lie close to an extended center manifold, on which the dynamics reduce to a driftless scalar problem solved in closed form by a scalar Riccati equation. Under weak control, the reduced problem on the same manifold is fully nonlinear, and the unfolding of the cusp enters only through the constant and linear coefficients of the reduced drift. In the last two regimes, the reduced problem depends on the system only through $\sigma$, $\mu$, $\nu$ and $\beta_c$, each a finite combination of derivatives of $r^i$ at the organizing center.

In the weak regime, a canonical change of variables reduces the extremal flow to motion in a potential. On the pitchfork ray, the bifurcation diagram of this flow organizes the possible extremals, and its saddle-node curve at large $\kappa$ marks the transition to the intermediate regime. In the long-horizon limit, extremals are sequences of heteroclinic layers whose positions follow from a far-field interaction law with explicit coefficients throughout the bistable part of the unfolding. The resulting profiles agree with numerically optimized solutions for both testbeds.

The application of the reduction to finite-horizon extremals and the boundary-layer laws are leading-order and formal. The second-order results (Theorems~\ref{thm:second_order} and~\ref{thm:local_optimality}, Corollary~\ref{cor:turning}) are rigorous, conditional on the layered structure (H) of the extremal, which we establish only formally. Under (H), the second variation of the full problem is coercive on $L^2$, uniformly in $\varepsilon$, in the strong and intermediate regimes. In the weak regime it is coercive whenever the Jacobi field $y$ has no zero on $[0,T_\tau]$. In each of these cases the extremal is a strict weak local minimizer for small $\varepsilon$. A zero of $y$ in $(0,T_\tau)$ rules out local optimality. For the reduced problems these statements are exact. Selecting the global optimum among competing locally optimal extremals still requires comparing their costs. The analysis assumes that the control couples directly to the critical mode; the complementary case falls under the control-bifurcation framework of Kang and Krener \cite{kang1998bifurcation1,krener2004control}. The numerical validation is restricted to the pitchfork ray.

The approach is not specific to the cusp: the reduction in Theorem~\ref{thm:reduction} rests on the center-manifold and Hamiltonian structure of the optimality system, and the cusp enters only through the distinguished scaling. Table~\ref{tab:other_bifs} lists this scaling for other structurally stable local bifurcations. For the cuspoid catastrophes $A_k$ \cite{poston1978catastrophe}, with reduced drift $f(A) = \sum_{j=0}^{k-2} \sigma_j h_j A^j + \nu A^k$, the weights $A \sim \varepsilon^{1/2}$ and $h_j \sim \varepsilon^{(k-j)/2}$ place the weak regime at $\omega = k-1$ and the driftless intermediate regime at $0 < \omega < k-1$, and leave the strong regime unchanged; for the fold ($k=2$), the intermediate regime falls outside the half-power expansions used here. The transcritical and symmetric pitchfork bifurcations are structurally stable only within the classes of systems with a persistent equilibrium or a reflection symmetry \cite{golubitsky1985singularities}; there the transcritical has the scaling of the fold, and the pitchfork that of the cusp on the ray $\tilde h_1 = 0$. Higher-dimensional center subspaces, as at the umbilic catastrophes and at Hopf points, would require relaxing Assumptions~\ref{assump:bifurcation} and~\ref{assump:diag}.

\begin{table}[htbp]
    \centering
    \small
    \setlength{\tabcolsep}{3pt}
    \begin{tabular}{lcllccc}
        \toprule
        & & & & \multicolumn{2}{c}{Weak regime} & Intermediate \\
        \cmidrule(lr){5-6}
        Bifurcation & Codim. & Reduced drift & Weights & $\omega$; $u$ & $T$ & $\omega$ \\
        \midrule
        Fold ($A_2$) & 1 & $\sigma h + \nu A^2$ & $h \sim \varepsilon$ & $1$; $\varepsilon$ & $\varepsilon^{-1/2}$ & -- \\
        Swallowtail ($A_4$) & 3 & $\sum_{j=0}^{2} \sigma_j h_j A^j + \nu A^4$ & $h_j \sim \varepsilon^{(4-j)/2}$ & $3$; $\varepsilon^{2}$ & $\varepsilon^{-3/2}$ & $1, 2$ \\
        Hopf$^\dagger$ & 1 & $(i\Omega + \lambda h)\, z + \ell\, |z|^2 z$ & $h \sim \varepsilon$ & $2$; $\varepsilon^{3/2}$ & $\varepsilon^{-1}$ & $1$ \\
        \bottomrule
    \end{tabular}
    \caption{Distinguished scalings for structurally stable local bifurcations, with state amplitude $O(\varepsilon^{1/2})$ in each case. The penalties scale as $\chi = \varepsilon^{\omega/2}\tilde\chi$ and $\zeta = \varepsilon^{\omega}\tilde\zeta$ (Eq.~\eqref{eq:alpha_rescaling}), the control as $u = O(\varepsilon^{(\omega+1)/2})$, and the strong regime $\omega = 0$ is the same in every case. Intermediate regimes are listed for integer $\omega$, i.e.\ within the half-power expansions of Sec.~\ref{sec:rescaled_system}. $^\dagger$Complex amplitude $z$ of the critical pair $\pm i\Omega$, with $\operatorname{Re}\ell \neq 0$; for $\omega \geq 1$ the reduced problem is obtained after averaging over the fast phase $\Omega t$. This entry is formal power counting.}
    \label{tab:other_bifs}
\end{table}

Two further extensions are natural. First, the connection to large-deviation theory noted in Sec.~\ref{sec:boundary_layer} suggests that weak noise can be incorporated within the same reduced Hamiltonian, although noise also modifies the bifurcation of the unforced dynamics \cite{crauel1998additive}. Second, formulations that penalize sensitivity to uncertainty in the unfolding parameters can be posed within the same framework.

\appendix

\section{Expansion details}
\label{app:expansions}

The coefficients of Sec.~\ref{sec:organizing_center} are, for $a \in \{1,2\}$,
\begin{align}
    \sigma_a &= w_c^m r^m{}_{h_a}, \nonumber\\
    \mu_a &= w_c^m r^m{}_{h_a i} v_c^i - \sum_s \frac{( \overline{w_s}^m r^m{}_{h_a} )(w_c^n r^n{}_{ij} v^i_c v_s^j)}{\lambda_s}, \label{eq:mu_a}\\
    \nu &= \frac{w_c^m r^m{}_{ijk} v^i_c v_c^j v_c^k}{6}  - \sum_s \frac{( \overline{w_s}^m r^m{}_{ij} v^i_c v_c^j)(w_c^n r^n{}_{kl} v^k_c v_s^l)}{2\lambda_s}.
     \label{eq:nu}
\end{align}
In $\mu_a$, the first term is the direct contribution. The sum is the indirect contribution: $h_a$ displaces the stable modes by $-(\overline{w_s}^m r^m{}_{h_a})/\lambda_s$, and that displacement couples back to the critical mode through the quadratic term. In $\nu$, the sum is the corresponding correction from stable modes slaved to $A^2$.

With the expansions \eqref{eq:c_expansion}--\eqref{eq:u_expansion}, the right-hand side expands as
\begin{align}
    r^i(c^j, u^\alpha, h) &= \varepsilon^{1/2}\, \left[ r^i{}_j\, c^{(1), j} + r^i{}_\alpha u^{(1), \alpha} \right] \nonumber\\
    &+ \varepsilon \Big[ r^i{}_j\, c^{(2), j} + r^i{}_\alpha u^{(2), \alpha} + \frac{1}{2}\, r^i{}_{jk}\, c^{(1), j} c^{(1), k} + r^i{}_{h_2}\tilde h_2 \nonumber\\
    &\qquad + r^i{}_{\alpha j} u^{(1), \alpha} c^{(1), j}  + \frac{1}{2} r^i{}_{\alpha \beta}  u^{(1), \alpha} u^{(1), \beta}  \Big] \nonumber\\
    &+ \varepsilon^{3/2} \bigg[ r^i{}_j\, c^{(3), j}  +  r^i{}_{jk}\, c^{(1), j} c^{(2), k}  +  \frac{1}{6}\, r^i{}_{jkl}\, c^{(1), j} c^{(1), k} c^{(1), l} + r^i{}_{h_1}\tilde h_1 \nonumber  \\
    &+  r^i{}_{j h_2}\, c^{(1), j}\tilde h_2 + r^i{}_\alpha u^{(3), \alpha}   + r^i{}_{\alpha j} u^{(2), \alpha} c^{(1), j} + r^i{}_{\alpha k} u^{(1), \alpha} c^{(2), k} \nonumber \\
    &+  r^i{}_{\alpha \beta}  u^{(2), \alpha} u^{(1), \beta} +  r^i{}_{\alpha h_2} u^{(1), \alpha}\tilde h_2 \nonumber \\
    &  + \frac{1}{6}r^i{}_{\alpha \beta \gamma}  u^{(1), \alpha} u^{(1), \beta} u^{(1), \gamma}  + \frac{1}{2}r^i{}_{\alpha \beta j}  u^{(1), \alpha} u^{(1), \beta}  c^{(1), j}  + \frac{1}{2}r^i{}_{jk\alpha}  c^{(1), j} c^{(1), k} u^{(1), \alpha} \bigg] + O(\varepsilon^2).
    \label{eq:r_expansion}
\end{align}
Each derivative is a partial derivative of $r^i$ with respect to the state components $c^j$, the control components $u^\alpha$, or the parameters $h_a$, evaluated at the organizing center.

The stable-mode amplitudes below involve the critical-to-stable couplings
\begin{equation*}
\rho_s = \overline{w_s}^i r^i{}_{h_2}, \qquad
\eta_{sc} = \overline{w_s}^m r^m{}_{jk}\, v_{c}^j v_c^k, \qquad
\eta_{cs} = w_c^m r^m{}_{ij}\, v_s^i v_c^j, \qquad
\beta_{s} = \overline{w_s}^m r^m{}_\alpha r^n{}_\alpha w_c^n ,
\end{equation*}
and $A^*_s = v_s^i c^{*(1),i}$ is the stable-mode projection of the target.

\subsection{Intermediate regime, $\omega = 1$}
\label{app:intermediate}

We take $\omega = 1$, so that $\chi = \varepsilon^{1/2}\tilde\chi$, $\zeta = \varepsilon\tilde\zeta$, $p^i = \varepsilon^{1/2}q^i$ and $\tau = \varepsilon^{1/2}t$. Substituting the expansions \eqref{eq:c_expansion}--\eqref{eq:u_expansion} into \eqref{eq:rescaled_state}--\eqref{eq:rescaled_control} and collecting powers of $\varepsilon$ gives the following hierarchy.

At lowest order, $O(\varepsilon^{1/2})$, time derivatives drop out:
\begin{align}
r^i{}_j c^{(1), j} + r^i{}_\alpha u^{(1), \alpha} &= 0 ,
\label{eq:c1_pde} \\
- r^j{}_i q^{(1), j} & = 0,
\label{eq:q1_pde} \\
u^{(1), \alpha} &= 0.
\end{align}
Since $u^{(1)} = 0$, state and costate decouple, and the lowest-order solutions lie in the null spaces of the Jacobian and of its adjoint, respectively:
\begin{equation}
    c^{(1),i} = A(\tau) v_c^i, \qquad
    q^{(1),i} = Q(\tau) w_c^i,
\end{equation}
where the amplitudes $A$ and $Q$ are determined by solvability at the next order.

At $O(\varepsilon)$, the hierarchy gives
\begin{align}
\frac{d c^{(1),i}}{d \tau} &=   r^i{}_j\, c^{(2), j}  + \frac{1}{2}\, r^i{}_{jk}\, c^{(1), j} c^{(1), k} + r^i{}_{h_2}\tilde h_2  + r^i{}_\alpha u^{(2), \alpha},
\label{eq:state_order1} \\
\frac{d q^{(1), i}}{d \tau} &= -\tilde{\zeta}(c^{(1), i} - c^{*(1), i}) - r^j{}_i q^{(2), j} - r^j{}_{ik}\,  c^{(1), k} q^{(1), j},
\label{eq:costate_order1} \\
u^{(2), \alpha} &= - r^i{}_\alpha q^{(1), i},
\end{align}
where terms containing $u^{(1)} = 0$ have been dropped. By the Fredholm alternative, projecting the state equation onto the left critical eigenmode $w_c^i$ and the costate equation onto $v_c^i$ gives the solvability conditions \eqref{eq:2nd_opt}. Here $w_c^i v_c^i = 1$ gives $A'$ and $-\beta_c Q$ in the state equation, and $\|v_c\| = 1$ gives $v_c^i(c^{(1),i} - c^{*(1),i}) = A - A^*$ in the costate equation. Three further projections vanish: $w_c^i r^i{}_{jk}\, v_{c}^j v_c^k$ in the state equation and the same coefficient $v_c^i w_c^j r^j{}_{ik} v_c^k$ in the costate equation, both by Assumption~\ref{assump:cusp}, and $w_c^i r^i{}_{h_2}$, by the parameter coordinates \eqref{eq:normalized_unfolding}. The constant $\sigma\tilde h_1$ first appears at $O(\varepsilon^{3/2})$.

For the stable-mode components, write $c^{(2),i} = \sum_\ell a_\ell v_\ell^i$ and $q^{(2),i} = \sum_\ell b_\ell\, \overline{w_\ell}^{\,i}$. Projecting \eqref{eq:state_order1} onto $\overline{w_s}^{\,i}$ and \eqref{eq:costate_order1} onto $v_s^i$ gives
\begin{align}
 a_s(\tau) &=  \frac{1}{\lambda_s} \left[ -\frac{A(\tau)^2}{2} \eta_{sc} - \rho_s \tilde h_2 + \beta_{s} Q(\tau) \right],
\label{eq:a_s}\\
    b_s(\tau) &= \frac{1}{\lambda_s} \left[ \tilde{\zeta}\left(A^*_s - A(\tau)\,v_s^i v_c^i\right) - \eta_{cs} A(\tau) Q(\tau)  \right].
    \label{eq:b_s_int}
\end{align}
Higher orders follow in the same way; for constant $c^{*i}$ they reduce to combinations of exponentials.

\subsection{Weak regime, $\omega = 2$}
\label{app:weak}

We take $\omega = 2$, so that $\chi = \varepsilon\tilde\chi$, $\zeta = \varepsilon^2\tilde\zeta$, $p^i = \varepsilon q^i$ and $\tau = \varepsilon t$. The rescaled optimality conditions \eqref{eq:rescaled_state}--\eqref{eq:rescaled_control} read
\begin{align}
    \varepsilon \frac{d c^i}{d\tau} &= r^i(c^j, u^\alpha, h), \\
    \varepsilon \frac{d q^i}{d\tau} &= - \varepsilon  \tilde{\zeta} (c^i - c^{*i}) - r^j{}_i q^j, \\
    0 &= u^\alpha + \varepsilon\, r^i{}_\alpha q^i.
\end{align}

As in the intermediate regime, the lowest order, $O(\varepsilon^{1/2})$, gives
\begin{equation}
    r^i{}_j c^{(1), j} = 0, \qquad
    - r^j{}_i q^{(1), j} = 0, \qquad
    u^{(1), \alpha} = 0,
\end{equation}
so $c^{(1),i} = A(\tau) v_c^i$ and $q^{(1),i} = Q(\tau) w_c^i$. The amplitudes $A$ and $Q$ are again fixed by solvability, which now enters at $O(\varepsilon^{3/2})$.

At the next order, $O(\varepsilon)$, the control is again absent:
\begin{align}
  r^i{}_j\, c^{(2), j}  + \frac{1}{2}\, r^i{}_{jk}\, c^{(1), j} c^{(1), k} + r^i{}_{h_2}\tilde h_2  &= 0,
\label{eq:weak_state_order2} \\
  r^m{}_i q^{(2), m} + r^m{}_{ij}\,  c^{(1), j} q^{(1), m} &= 0,
\label{eq:weak_costate_order1} \\
u^{(2), \alpha} &= 0.
\end{align}
The critical projections of \eqref{eq:weak_state_order2} and \eqref{eq:weak_costate_order1} vanish identically by Assumption~\ref{assump:cusp} and \eqref{eq:normalized_unfolding}, so these equations are solvable. Writing $c^{(2),i} = a(\tau) v_c^i + \sum_{s} a_s(\tau) v_s^i$ and $q^{(2),i} = b(\tau) w_c^i + \sum_{s} b_s(\tau)\, \overline{w_s}^{\,i}$, and projecting onto $\overline{w_s}^{\,i}$ and $v_s^i$, the stable-mode amplitudes are
\begin{equation}
    a_s(\tau) = - \frac{1}{\lambda_s} \left[ \frac{1}{2 } \eta_{sc} A(\tau)^2 + \rho_s\tilde h_2 \right], \qquad
    b_s(\tau) = - \frac{1}{\lambda_s}\, \eta_{cs} A(\tau) Q(\tau).
    \label{eq:as_weak}
\end{equation}
The critical contributions $a(\tau)$ and $b(\tau)$ are undetermined at this stage, but they do not enter the amplitude equations at leading order.

At $O(\varepsilon^{3/2})$, omitting terms that vanish by the lower-order results, we have
\begin{align}
    \frac{d c^{(1),i}}{d \tau} &=  r^i{}_j\, c^{(3), j}  +  r^i{}_{jk}\, c^{(1), j} c^{(2), k}  +  \frac{1}{6}\, r^i{}_{jkl}\, c^{(1), j} c^{(1), k} c^{(1), l} + r^i{}_{h_1}\tilde h_1 +  r^i{}_{j h_2}\, c^{(1), j}\tilde h_2  +  r^i{}_\alpha u^{(3), \alpha},    \\
    \frac{d q^{(1), i}}{d \tau}     &= \tilde{\zeta}( c^{*(1), i} - c^{(1), i} )- r^m{}_i q^{(3), m} - r^m{}_{ik}\, c^{(1), k} q^{(2), m} - r^m{}_{ik}\, c^{(2), k} q^{(1), m} \nonumber \\
    &\quad -\frac{1}{2}r^m{}_{ijk} c^{(1), j}  c^{(1), k} q^{(1), m} -  r^m{}_{i h_2}\tilde h_2\, q^{(1), m}, \\
    u^{(3), \alpha} &= - r^i{}_\alpha q^{(1), i}.
\end{align}
Projecting onto the critical mode, and substituting \eqref{eq:as_weak} for the stable-mode components of $c^{(2)}$ and $q^{(2)}$, yields the solvability conditions \eqref{eq:Aweak}--\eqref{eq:Qweak}. This substitution generates the stable-mode sums in $\mu = \mu_2$ and $\nu$, Eqs.~\eqref{eq:mu_a}--\eqref{eq:nu}. The constant $\sigma\tilde h_1$ arises from the projection $w_c^i r^i{}_{h_1}\tilde h_1$ in the state equation. Its counterpart in the costate equation, $r^m{}_{i h_1} h_1 q^m$, is $O(\varepsilon^2)$ and does not appear. This is consistent with the Hamiltonian structure of Theorem~\ref{thm:reduction}: a constant drift term has no derivative with respect to $A$.

\section{Numerical details}
\label{app:numerics}

All code used for the computations and figures in this paper is provided as supplementary material.

Both testbeds below are invariant under a reflection that reverses the critical mode. They therefore lie on the ray $\tilde h_1 = 0$ of the unfolding, with a shifted bifurcation parameter of each model playing the role of $h_2$, as specified below.

In Figs.~\ref{fig:lqr}--\ref{fig:nonlinear}, the asymptotic solution is the trajectory of the full system driven by the leading-order control law of the corresponding regime, and the reported error is the full-state $L_\infty$ difference from the numerical optimum. Throughout, a numerical optimum is a stationary point of the discretized problem found by local optimization from an asymptotic warm start; its minimality is not certified. In the weak regime, the bistable switch is driven by the Hamiltonian level-set feedback law of Sec.~\ref{sec:weak_regime}, and the Turing model by the solution of the reduced boundary-value problem. Since the level $E$ is fixed only implicitly by the boundary conditions, for the bistable switch we take $E$ as the median of $H_{\text{reduced}}$ along the reduced amplitudes $(A,Q)$ of the numerical optimum, and integrate the feedback law at this level.

\subsection{Bistable switch}

We consider the classic systems-biology motif of a bistable switch, using a mutually inhibitory circuit adapted from Ref.~\cite{ferrell2012bistability} (inspired by Delta--Notch signaling). The controlled dynamics are
\begin{align}
  \frac{dx_1}{dt} &= \frac{k}{1 + (I x_2)^n} - d\bigl(1 + U_1(t)\bigr) x_1, \\
  \frac{dx_2}{dt} &= \frac{k}{1 + (I x_1)^n} - d\bigl(1 + U_2(t)\bigr) x_2,
\end{align}
with production rate $k=2$, degradation rate $d=2$, and Hill coefficient $n=4$. The controls $\bigl(U_1(t),U_2(t)\bigr)$ modulate the effective degradation rates of the two species. For $U_1=U_2=0$, the system undergoes a symmetry-breaking pitchfork bifurcation as $I$ varies, and we take $h_2 = I - I_c$, i.e.\ $I = I_c + \varepsilon$ ($\tilde h_2 = 1$). The critical value is
\begin{equation}
  I_c = \frac{4}{3\cdot 3^{1/4}} \approx 1.01311,
\end{equation}
with the symmetric steady state at criticality $x_c=3/4$. A weakly nonlinear analysis yields the normal-form coefficients $\mu\approx 2.96117$, $\nu\approx -2.22222$, and the critical-mode control coupling $\beta_c=9/4$.

We consider the problem of driving the system from one asymmetric stable equilibrium to the other, in all three control regimes and over a range of $\varepsilon$. In the computations reported here, the state-tracking coefficient is $\zeta$ (scaled by regime as described above), and the terminal coefficient $\chi$ is held at zero. The boundary-layer comparison of Fig.~\ref{fig:boundary_layer_comparison} instead uses $\chi = 1000\,\varepsilon$ ($\tilde\chi = 1000$), which gives the values of $\hat\chi = \beta_c\tilde\chi/|\hat\mu|$ quoted there. The time horizon and step count are set by regime: $T=50$, $N_t=10^4$ (strong); $T=100$, $N_t=10^4$ (intermediate); and $T=2000$, $N_t=2\times10^4$ (weak). For both testbeds $T$ is fixed while $\varepsilon$ varies, so $T_\tau = \varepsilon^{\omega/2}T$ changes along each sweep, and the error rates in Figs.~\ref{fig:lqr}--\ref{fig:nonlinear} are empirical for these families. The boundary-layer comparison uses $\varepsilon = 10^{-2}$ and $T = 2000$ for the bistable switch and $T = 750$ for the Turing model, with $N_t = 4\times10^4$. It is computed not by the direct method below but by solving the full-system Pontryagin boundary-value problem with the implicit midpoint rule, with damped Newton iteration on the state and costate at every node. This scheme is the discrete optimality system of the midpoint-rule discretization of the optimal control problem, so the computed extremal is a stationary point of the discrete cost. The scheme is symplectic, so for small steps and accurate Newton solves $H$ is nearly conserved, to $O(\Delta t^2)$, on bounded trajectories. The results in Figs.~\ref{fig:lqr}, \ref{fig:feedback}, and \ref{fig:nonlinear} show that the asymptotic laws capture the leading-order structure of the numerically optimized trajectories in each regime.

The discretized optimal control solution is computed using a Python implementation with the JAX automatic-differentiation framework. The time grid is $t_k=k\Delta t$ for $k=0,\dots,N_t$, with $\Delta t=T/N_t$. For dynamics $\dot{x}=r(x,u)$, we advance the state using first-order implicit Euler,
\begin{equation}
  x_{k+1} = x_k + \Delta t\, r(x_{k+1},u_k),
\end{equation}
which we evaluate via a fixed-point iteration initialized at $x_k$,
\begin{equation}
  x_{k+1}^{(m+1)} = x_k + \Delta t\, r\bigl(x_{k+1}^{(m)},u_k\bigr),
\end{equation}
using five iterations per time step, both in the objective/gradient forward solve and when integrating the reported optimal and feedback-driven trajectories. The truncated iteration defines a fixed discrete map, not an exact implicit Euler solve, and automatic differentiation is exact for this map. The discretized objective is
\begin{equation}
  J_h(U) = \frac{1}{\varepsilon}\left[\frac{\Delta t}{2}\sum_{k=0}^{N_t-1}\left(\zeta\,\|x_{k+1}-x_{k+1}^*\|_2^2 + \|u_k\|_2^2\right) + \frac{\chi}{2} \|x_{N_t}-x_{N_t}^*\|_2^2\right].
\end{equation}
We compute gradients by reverse-mode automatic differentiation through the unrolled iteration in JAX, avoiding hand-derived adjoint discretizations. Optimization proceeds in two stages: (i) gradient descent with a warm start from the asymptotic approximation, followed by (ii) refinement via a Newton--Krylov solve. The overall factor $1/\varepsilon$ does not change the minimizer; it is included to improve the scaling of the gradient descent.

\subsection{Reaction--diffusion model}

As a second, structurally distinct testbed, we consider a reduced four-cell activator--substrate reaction--diffusion system of Schnakenberg type on a one-dimensional domain with Neumann (no-flux) boundary conditions. Let $u\in\mathbb{R}^4$ denote the activator and $v\in\mathbb{R}^4$ the substrate. The substrate here is slaved to the activator to keep the example simple:
\begin{align}
  \frac{du}{dt} &= \mat{L}u + \gamma\bigl(u \odot u \odot v - u\bigr) + \mathbf{e}_1\,U(t), \\
  0 &= \delta_{\mathrm{sub}}\,\mat{L}v + \mathbf{1} - u \odot u \odot v,
\end{align}
where $\odot$ denotes elementwise multiplication, $\mathbf{1}=(1,1,1,1)^\top$, and $\mat{L}$ is the $4\times4$ Neumann finite-volume Laplacian on a grid of spacing $\Delta x=0.25$. The scalar control $U(t)$ acts through $\mathbf{e}_1$, i.e.\ a localized source at the first grid cell. The substrate diffusivity ratio is fixed at $\delta_{\mathrm{sub}}=0.6$, and the reaction rate $\gamma$ plays the role of the bifurcation parameter (analogous to $I$ above).

For $U=0$, the homogeneous steady state $u^*=\mathbf{1}$ is always an equilibrium. It undergoes a symmetry-breaking (Turing) instability as $\gamma$ increases past a mode-dependent threshold. Writing $\lambda^{\mathsf{L}}_k = 32\bigl(\cos(k\pi/4)-1\bigr)$ for the eigenvalues of $\mat{L}$ restricted to the $k$th Neumann cosine mode, the threshold for the first nontrivial mode ($k=1$) is
\begin{equation}
  \gamma_c = |\lambda^{\mathsf{L}}_1|\,\frac{1+\delta_{\mathrm{sub}}|\lambda^{\mathsf{L}}_1|}{\delta_{\mathrm{sub}}|\lambda^{\mathsf{L}}_1|-1} \approx 13.4269,
\end{equation}
with associated critical mode shape $v_c^i \propto \cos\!\bigl(\pi(i-\tfrac12)/4\bigr)$, normalized so that $\|v_c\|_2=1$. A weakly nonlinear analysis of the reduced four-cell model yields, with $h_2 = \gamma/\gamma_c - 1 = \varepsilon$ ($\tilde h_2 = 1$), $\mu\approx 9.37258$, $\nu\approx -3.17077$, and the critical-mode control coupling $\beta_c=(v_c^1)^2\approx 0.42678$, where $v_c^1$ is the first (control-actuated) component of the critical mode. As in the bistable-switch example, $\mu\nu<0$, consistent with a supercritical pitchfork; in the notation of Sec.~\ref{sec:nonlinear_analysis}, this is the ray $b = 0$ with $a = -1$.

We pose the same optimal control problem as for the bistable switch, at fixed $\gamma = \gamma_c(1+\varepsilon)$: drive the system from one patterned branch to the other, the two branches being obtained by following the critical mode to either side of $u^*$. The numerical scheme reflects the semilinear structure of the reduced activator--substrate system. Diffusion is treated implicitly and the (already slaved) reaction term explicitly (an IMEX splitting), which avoids the fixed-point iteration used above. At each time step $k$, we eliminate the slaved substrate by solving
\begin{equation}
  \bigl(\mathrm{diag}(u_k \odot u_k) - \delta_{\mathrm{sub}}\,\mat{L}\bigr)\,v_k = \mathbf{1},
\end{equation}
then advance the activator via
\begin{equation}
  \bigl(\mat{I} - \Delta t\,\mat{L}\bigr)\,u_{k+1} = u_k + \Delta t\left[\gamma\bigl(u_k \odot u_k \odot v_k - u_k\bigr) + \mathbf{e}_1 U_k\right].
\end{equation}
The discretized objective has the same form as $J_h(U)$ above, with $x_k \to u_k$ and $u_k\to U_k$. Gradients are again obtained by reverse-mode automatic differentiation through the IMEX unroll, and the same two-stage (gradient descent, then Newton--Krylov) minimization is used, warm-started from the corresponding asymptotic solution. The time horizon and step count are set by regime: $T=100$, $N_t=5\times 10^3$ (strong/LQR); $T=200$, $N_t=2\times 10^4$ (intermediate/feedback); and $T=500$, $N_t=3\times 10^4$ (weak/manifold-BVP warm start).

\subsection{Bifurcation analysis}
\label{app:bifurcation_numerics}

The bifurcation sets of Figs.~\ref{fig:local_bifs} and~\ref{fig:codim2_tracks} are computed in Mathematica, and the code is included in the supplementary material. Write $W = g^2/2$, so that $V_X = \kappa(X - X^*) + W'(X)$ and $V_{XX} = \kappa + W''(X)$. The saddle-node set is then parametrized by its degenerate point $X$, with $\kappa = -W''(X)$ and $X^* = X - W'(X)/W''(X)$. For the heteroclinic set, two critical points $x \neq y$ of a common potential fix $\kappa$ and $X^*$ through divided differences of $W'$, and the equal-level condition $V(x) = V(y)$ becomes a polynomial in $(x, y)$. Its zero set is traced by contouring and kept where $x$ and $y$ are both minima of $V$ and every critical point between them lies above their common level. The codimension-two points at $b = 0$ are solved directly as polynomial systems. To continue them in $b$, we use that each defining condition is linear in $b$ to eliminate it, and for the heteroclinic--saddle-node intersections we eliminate the third critical point by a resultant. This yields either an explicit parametrization (the cusps) or a planar algebraic curve in the critical points $(x, y)$. These curves are traced by contouring and mapped to $(X^*, b, \kappa)$, so turning points in $b$ require no branch tracking. The codimension-three points that end the curves and the $\kappa \to 0$ asymptotes of the double saddle-node curve are obtained in closed form. The phase portraits of Fig.~\ref{fig:local_bifs} are streamlines of the canonical flow at $b = 0$.

\bibliography{references}

\section*{Competing Interests}
There are no competing interests to declare.

\section*{Funding}
This work was funded through internal grants from UC San Diego.

\section*{Data Availability}
No empirical data sets were involved in this study. All code used for the simulations, the bifurcation analysis and the figures is provided as supplementary material.

\section*{AI Declaration}
Generative AI tools were used in this work under the author's direction and oversight. GitHub Copilot (in VS Code) assisted in developing and debugging the research code for the simulations, bifurcation analysis and figures, and in \TeX{} formatting. The code was reviewed by the author and validated against independent checks, including closed-form values of the critical parameters and normal-form coefficients of both testbeds. Claude (Anthropic) was used for literature search, for proofreading and language editing, and for drafting revised passages from the author's instructions and referee comments; every mathematical statement in these passages was checked by the author. All bibliographic references were verified against publisher records. No figure was generated by AI. The author takes full responsibility for all research decisions, results and conclusions in this manuscript.

\section*{Acknowledgments}

We thank Vishal Patil, Natalia Komarova, and Martin M\"onnigmann for helpful discussions.

\end{document}